\documentclass[11pt]{amsart}
\usepackage{epsfig, amssymb}

\newtheorem{theorem}{Theorem}[section]
\newtheorem{proposition}[theorem]{Proposition}
\newtheorem{lemma}[theorem]{Lemma}
\newtheorem{corollary}[theorem]{Corollary}

\theoremstyle{definition}
\newtheorem{definition}[theorem]{Definition}

\theoremstyle{remark}

\numberwithin{equation}{section}

\newcommand{\be}{\begin{equation}}
\newcommand{\ee}{\end{equation}}

\newcommand{\zbar}{\overline{z}}
\newcommand{\h}{\hbar}

\newcommand{\N}{{(N)}}

\newcommand{\sgn}{\mbox{\small sgn}}

\newcommand{\dint}{\int\hspace{-8pt}\int}

\newcommand{\ket}[1]{|#1{\rangle}}

\newcommand{\wf}{\mbox{WF}}
\newcommand{\bbC}{{\mathbb C}}

\newcommand{\bbZ}{{\mathbb Z}}
\newcommand{\bbR}{{\mathbb R}}
\newcommand{\bbP}{{\mathbb P}}
\newcommand{\calA}{{\mathcal A}}
\newcommand{\calB}{{\mathcal B}}

\newcommand{\calL}{{\mathcal L}}
\newcommand{\calM}{{\mathcal M}}

\newcommand{\calH}{{\mathcal H}}

\newcommand{\calF}{{\mathcal F}}
\newcommand{\calT}{{\mathcal T}}

\newcommand{\calR}{{\mathcal R}}
\newcommand{\calS}{{\mathcal S}}

\newcommand{\calO}{{\mathcal O}}

\newcommand{\calZ}{{\mathcal Z}}

\DeclareMathOperator{\tr}{Tr}

\DeclareMathOperator{\Tr}{tr}

\newcommand{\inner}[2]{\langle#1,#2\rangle}
\newcommand{\norm}[1]{\lVert#1\rVert}

\newcommand{\PB}[2]{\{  #1\,,\,#2 \}}

\newcommand{\del}{\partial}

\newcommand{\hffm}{{\textstyle \bigwedge^{1/2}}}

\newcommand{\nchffm}{{\textstyle \overline{\bigwedge}^{-1/2}}}

\begin{document}

\title{On the pseudospectra of Berezin-Toeplitz operators}
\author{D. Borthwick}
\address{Department of Mathematics and Computer Science\\
Emory University\\Atlanta, GA 30322}
\email{davidb@math.emory.edu}
\thanks{Supported in part by NSF grant DMS-0204985.}
\author{A. Uribe}
\address{Mathematics Department\\
University of Michigan\\Ann Arbor, Michigan 48109}
\email{uribe@umich.edu}
\thanks{Supported in part by NSF grant DMS-0070690.}

\date{14 April 2003}

\maketitle

\tableofcontents

\section{Introduction}

The study of pseudospectra of operators is a very active area of research,
of importance in applied mathematics and numerical analysis. 
We recall that if $\epsilon > 0$,
the $\epsilon$-pseudospectrum of an operator $Q$ acting on a Hilbert space
is the set of complex numbers, $\lambda\in\bbC$, such that
\[
\norm{(Q-\lambda I)^{-1}} \geq \frac{1}{\epsilon},
\] 
or, equivalently, the set of $\lambda$ such that
\begin{equation}\label{1a}
\inf_{\psi\not= 0}\,\frac{\norm{(Q-\lambda I)(\psi)}}{\norm{\psi}} \leq
\epsilon.
\end{equation}
The importance of this set in applied settings stems from the following two
facts: 
\begin{enumerate}
\item If $\epsilon$ is very small it may be difficult to distinguish
the $\epsilon$-pseudospectrum from the spectrum of $Q$, 
\item If $Q$ is strongly non normal the $\epsilon$-pseudospectrum generally
is much bigger than the spectrum, even if $\epsilon$ is very small.
\end{enumerate}
In case $Q$ is a differential or, more generally, a pseudodifferential
operator with a small parameter, $\h$ (such as a Schr\"odinger operator),
it is natural to consider the
asymptotic behavior of the $\epsilon$-pseudospectrum of $Q$
where $\epsilon$ is related to $\h$, e.g.\ $\epsilon = O(\h^\infty)$. 
Dencker, Sj\"ostrand and Zworski have recently studied this problem
by microlocal techniques, \cite{DSZ}.  
Earlier results for Schr\"odinger operators were obtained
by E.\ B.\ Davies, \cite{Dav},
and  P.\ Redparth, \cite{Red}.

\medskip
On the other hand, in applications one often deals with large
{\em matrices}, and then it is natural to estimate the pseudospectrum
in terms of the size of the matrix.  In a recent paper, \cite{TC},
Trefethen and Chapman considered this problem for matrices
$T^\N = (T^N_{jl})$ where
\begin{equation}\label{1b}
T^\N_{jl} = f_{(l-j)\text{\small{mod}}N}(j/N),\quad
1\leq j,l\leq N.
\end{equation}
Here the $f_j$ are $1$-periodic coefficient functions . 
The main result of Trefethen and Chapman is that, under certain
assumptions including the following ``twist" condition,
\begin{equation}\label{1d}
\Im\Bigl(\frac{\partial f}{\partial p}/\frac{\partial f}{\partial x}\Bigr)(x_0,p_0)
< 0
\end{equation}
where $f(x,p) = \sum_j f_j(p) e^{ijx}$,
then $\lambda = f(x_0,p_0)$ is in the $\epsilon$-pseudospectrum of $T_N$
where $\epsilon = O(e^{-cN})$ for some $c>0$.

\medskip
The purpose of this note is to show that microlocal techniques can also
be applied to the study of the pseudospectra of matrices such as
(\ref{1b}) (and generalizations).  In this light we interpret the twist
condition (\ref{1d}) as H\"ormander's solvability condition
\begin{equation}\label{0}
\PB{\Re f}{\Im f}(x_0,p_0)<0
\end{equation}
on the Poisson bracket of the real and imaginary parts of
the symbol of a pseudodifferential operator.  Indeed if we define the
Poisson bracket of the variables $x$ and $p$ above to be one, it is
easy to check that (\ref{1d}) is precisely H\"ormander's condition.
The connection between H\"ormander's condition and pseudospectra was
first made by M.\ Zworski in \cite{Z}.  

In this paper we construct pseudomodes for Berezin-Toeplitz operators
under condition (\ref{0}) on the (smooth) symbol.  Our construction is
symbolic: in \S 2 we introduce spaces of Hermite distributions
containing the pseudomodes.
From the point of view of the symbolic calculus of these distributions,
condition (\ref{0}) is
exactly the condition on $\epsilon$ for the operator:
$\frac{d}{dx}+\epsilon x$ to have a kernel in the Schwartz space of
$\bbR$, namely $\epsilon > 0$.

Although we will discuss our results in detail in the next section, we
should mention some limitations of our work.  The methods of Trefethen
and Chapman apply to rough symbols, $f$, and they obtain exponentially
small error terms.  For analytic symbols, it is very likely that
exponentially small estimates (in the Toeplitz setting) can be achieved
by microlocal methods, as has been done in \cite{DSZ} for
pseudo-differential operators.  The problem of dealing with general
non-smooth symbols, however, is much more challenging.  Trefethen and
Chapman's main theorem includes a {\em global} condition on the symbol
(in addition to \ref{1b}), and they present compelling numerical
evidence that global conditions on non-smooth symbols are necessary for
the existence of ``good" pseudomodes (see \S 8 of \cite{TC}).  This is
a very interesting issue that we do not address here.  On the other
hand, our results for smooth symbols are fairly general and include a
number of cases not covered  by the results in \cite{TC} (e.\ g.\  the
``Scottish flag" matrix).  Furthermore, the pseudomodes we construct
are localized in phase space, sharpening the localization results of
\cite{TC}.

\medskip
We also mention that more straightforward microlocal methods can be
applied to the study of non-periodic versions of (\ref{1b}), along the
lines of the example in \S 4.1.  More generally, the Berezin-Toeplitz
operator calculus opens up the entire spectrum of phase-space methods
to study other problems associated with certain sequences of large
matrices, and we hope to return to some of these problems in the
future.

\subsection{The main results}

The general setting for B-T operators is a K\"ahler manifold, $X$, together
with a holomorphic hermitian line bundle, $L\to X$ whose curvature is the
symplectic form on $X$.  If $f:X\to\bbC$ is a ``classical Hamiltonian",
(a smooth function) consider the sequence of operators
$T_f=\{T_f^\N\,,\,N=1,2,\ldots\}$, acting on the space 
$ \calH_N $ of holomorphic sections of the tensor power $L^{\otimes N}$, 
defined by:
\[
\calH_N\ni\psi\mapsto \Pi_N(f\psi),
\]
where $\Pi_N: L^2(X,L^{\otimes N})\to\calH_N$ is orthogonal projection.
The sequence, $T_f$, is the primary example of a Berezin-Toeplitz operator.
More generally, one can allow $f$ to depend on $N$ as well, provided
the $N$-dependence admits an asymptotic expansion as $N\to\infty$
\[
f(x,N)\sim\sum_{j=0}^\infty N^{-j}\,f_j(x)
\]
in the $C^\infty$ topology.  The function $f_0$ is then called the
principal symbol of the operator.  We include explicit examples
of all this in \S 4.  General recent references for the theory of K\"ahler
quantization and Berezin-Toeplitz operators are \cite{Bor},
\cite{Ch}, \cite{Zel}.

For each $N$ $\calH_N$ is finite-dimensional, and for $N$ large, by the
Riemann-Roch theorem, $\dim \calH_N$ is a polynomial in $N$ of degree
one-half the dimension of $X$, $n:=\frac{1}{2}\text{dim }X$.  

The simplest cases are when $X$ is a either the torus or the complex
projective line, for which $n=1$.  Thus the parameter, $N$, is essentially
the dimension of $\calH_N$.  Moreover, for such $X$ the spaces
$\calH_N$ have a natural multiplicity-free representation of the
circle group, whose eigenvectors form a canonical basis of $\calH_N$.
We will write down explicitly the matrices corresponding to a B-T
operator on these spaces in the next section.  It turns out, for
example, that the matrices (\ref{1b}) are the matrices of B-T operators
on $X$ equal to the two-torus. 


Our main result is:

\begin{theorem}\label{I}
Let $T_f=\{T^\N\,,\,N=1,2,\ldots\}$ be a Berezin-Toeplitz operator
with smooth principal symbol $f:X\to\bbC$. 

\noindent
{\bf 0.} For all $\lambda\in\bbC$, 
\[
\inf_{\psi\in\calH_N}\,\frac{\norm{(T_f^\N-\lambda I)(\psi)}}{\norm{\psi}}
= \inf_{x\in X}\,|f(x)-\lambda| + O(1/\sqrt{N}).
\]

\noindent
{\bf 1.}   Assume that $\lambda = f(x_0)$
where $x_0\in X$ is such that
\begin{equation}\label{1f}
\PB{\Re f}{\Im f}(x_0) < 0.
\end{equation}
Then there exists a sequence of vectors $\{\psi_N\in\calH_N\}$ 
with microsupport precisely $\{x_0\}$ and such that
\begin{equation}\label{1g}
\frac{\norm{(T_f^\N-\lambda I)(\psi_N)}}{\norm{\psi_N}} = O(N^{-\infty}).
\end{equation}

\noindent
{\bf 2.} On the other hand, if $\lambda = f(x_0)$ and
\begin{equation}\label{1h}
\PB{\Re f}{\Im f}(x_0) > 0,
\end{equation}
then any sequence $\{\psi_N\in\calH_N\}$ such that (\ref{1g}) holds has
microsupport away from $\{x_0\}$.

\end{theorem}
We say a few words about the definition of microsupport in \S \ref{Two.one}.  

In the case of more than one degree of freedom
(i.\ e.\  if the dimension of $X$ is greater than two), there are
multiple pseudomodes under condition (\ref{1f}).  The level set
$f^{-1}(\lambda)$ is a symplectic manifold (at least near $x_0$) and 
we construct pseudomodes associated with any germ of isotropic submanifold
of it containing $x_0$.

\medskip
We will also prove an additional result, analogous to Theorem 4 in
\cite{DSZ}, whose hypotheses hold typically in case
$\lambda$ is on the boundary of the image of the principal symbol.
Let $\Re(f) = f_1$, $\Im(f) = f_2$, and for $I\subset \{1,2\}^m$ denote
by $f_I$ the repeated Poisson bracket:
\[
f_I = \Xi_{f_{i_1}}\,\Xi_{f_{i_2}}\,\ldots\,\Xi_{f_{i_{m-1}}} f_{i_m}
\]
where $\Xi_g$ denotes the Hamilton vector field of $g$.  Denoting the 
order
of the Poisson bracket by $|I| = m$, we define the order of a point 
$x\in X$ as
\begin{equation}
k(x) :=  \max\{j\in\bbZ\;;\;f_I(x)=0\ \text{for all}\ |I|\leq j\,\}.
\end{equation}

\begin{theorem}\label{II}
Let $\lambda\in\del\text{Image}(f)$ be such that:
\begin{enumerate}
\item $df_x\not= 0$ for every $x\in f^{-1}(\lambda)$.
\item The maximum order $k:=\max_{x\in f^{-1}(\lambda)} k(x)$ is finite.
\end{enumerate}
Then there exist $C,\ C_1>0$ such that
\begin{equation}\label{1i}
C_1N^{-1/2}\geq
\inf_{\psi\in\calH_N}\,\frac{\norm{(T_f^\N-\lambda I)(\psi)}}{\norm{\psi}}
\geq C N^{-\frac{k}{k+1}}.
\end{equation}
\end{theorem}

The first inequality in (\ref{1i}) follows from Part 0 of
Theorem \ref{I}; the second one follows from subelliptic estimates.
In the course of the proof of this Theorem we will also show
that, in general, the microsupport of a sequence of vectors
$\psi_N\in\calH_N$ minimizing the Rayleigh quotient in (\ref{1i})
(i.e.\ the optimal pseudomodes) is contained in $f^{-1}(\lambda)$, see
Proposition \ref{IIb}.  

The existence part of Theorem \ref{I} is proved by constructing
pseudomodes out of a class of distributions that will be defined in the
next section.  The proof of Theorems \ref{I} and \ref{II} appear in
\S 3, and \S 4 is devoted to examples.  We present some additional
results, including a description of the limit of the numerical range,
in \S 5.

\medskip\noindent {\sc Acknowledgments:}  We wish to thank Nick
Trefethen for sharing with us an early version of \cite{TC} and for 
encouraging remarks.  Many thanks also to Maciej Zworski and to Thierry Paul
for helpful conversations.

\section{Preliminaries}

\subsection{Setup and strategy}\label{Two.one}

Let $L\to X$ be as in the previous section, and 
let $P\subset L^*$ denote the unit circle bundle in the
dual of the line bundle $L$.  We denote by
\begin{equation}
L^2(P) = \bigoplus_{k\in\bbZ} L^2_k(P)
\end{equation}
the Fourier decomposition of functions on $P$ under the
action of the circle; explicitly $f\in L^2(P)$ is in 
$L^2_k(P)$ iff $f(e^{is}\cdot p) = e^{iks}f(p)$.  We will
also need the spaces
\begin{equation}
C^\infty_k(P) = L^2_k(P)\cap C^\infty(P).
\end{equation}
It is a tautology that $C^\infty_k(P)$ (resp. $L^2_k(P)$)
can be naturally identified with the space of
sections $C^\infty (X, L^{\otimes k})$ (resp. $L^2(X, L^{\otimes k})$.
We will henceforth identify these spaces without further comment.

$P$ is a strictly pseudoconvex domain, and under the natural action of the
circle group its Hardy space, $\calH$, splits into Fourier components
that are naturally isomorphic to the spaces of holomorphic sections
$\calH_N$:
\begin{equation}\label{wp0a}
\calH = \bigoplus_{N=0}^\infty \calH_N,\qquad \calH_N = H^0(X,L^{\otimes N}).
\end{equation}
We denote by $\Pi: L^2(P)\to\calH$ the Szeg\"o projector of $P$, and by
$\Pi_N: L^2(P)\to\calH_N$ the orthogonal projection onto the summand
$\calH_N$.
 
The precise structure of the singularities of $\Pi$ has been known
for some time, thanks to work of Boutet de Monvel and Sj\"ostrand.
We now recall the microlocal structure of $\Pi$, as described in \cite{BG}. 
Let $\alpha$ denote the
connection form on $P$, and let $\calZ\subset T^*P$ be the manifold
\[
\calZ = \{\,(p,r\alpha_p)\,;\,p\in P\ ,\ r>0\,\}.
\]
This is a conic symplectic submanifold of $T^*P$,
and $\Pi$ is a Fourier integral operator of Hermite type associated
with the canonical relation
\[
\calZ^\Delta := \{\,(\zeta,\zeta)\,;\,\zeta\in \calZ\,\}.
\]
We will say a few words below about the symbol of $\Pi$, referring 
to \cite{BG} for the general theory of Fourier integral operators
of Hermite type.  (See \cite{Zel} for a description of $\Pi$ as
a Fourier integral operator with complex phase.)

\medskip
The overall strategy of our proofs is the observation that
much of the asymptotic behavior of a sequence
$\{\psi_N\in C^\infty(X,L^{\otimes N})\}$ is encoded by the 
singularities of the distribution on $P$,
$\psi = \sum_{N=1}^\infty \psi_N\in C^{-\infty}(P)$.
For example, we have the following elementary result:

\begin{lemma}\label{Sobolev}
Given a sequence of vectors $\psi_N\in C^\infty(X,L^{\otimes N})$,
let $\psi = \sum_{N=1}^\infty \psi_N\in C^{-\infty}(P)$.
For each $s\in\bbR$, let $H_{(s)}(P)$ denote the Sobolev space
on $P$ consisting of distributions, $u$, such that 
$(\Delta_P+I)^{s/2}(u)\in L^2(P)$, with the norm
$\norm{u}_{(s)} = \norm{(\Delta_P+I)^{s/2}(u)}_{L^2}$.
Then the following are equivalent:
\begin{itemize}
\item[(a)] $ \psi\in C^\infty(P)$
\item[(b)] For all $s\in\bbR\ \norm{\psi_N}_{(s)} = O(N^{-\infty})$.
\end{itemize}
If, in addition, $\psi_N\in\calH_N$ for all $N$, then the above
are equivalent to:
\begin{itemize}
\item[(c)] $\norm{\psi_N}_{L^2} = O(N^{-\infty})$.
\end{itemize}
\end{lemma}
\begin{proof}
The spaces $C_k^\infty(P)$ are invariant under $\Delta_P$ and
orthogonal in $H_{(s)}$, so for each $s$ 
$ \norm{\psi}^2_{(s)} = \sum_N \norm{\psi_N}^2_{(s)}$.
It follows that (b) implies that $\psi\in H_{(s)}(P)$ for all $s$,
and therefore it implies (a) (and (c), of course).
Assuming (a) now, consider $D_\theta^k \psi$ where $k$ is a positive
integer.  Since this function is smooth it is in $H_{(s)}(P)$ for
each $s$, and therefore
$\norm{D_\theta^k \psi}_{(s)}^2 =
\sum_N N^{2k}\,\norm{\psi_N}_{(s)}^2 < \infty$,
which implies that $\norm{\psi_N}_{(s)} = O(N^{-k})$.  Therefore
(a) implies (b).  

Let us define an operator $\Delta_h$ by the identity
\begin{equation}\label{thm2ab}
\Delta_P = \Delta_h + D_\theta^2.
\end{equation}
Then $[\Delta_h,D_\theta]=0$ and the restriction of
$\Delta_h$ to $C^\infty_N(P)$ agrees exactly with the
Laplacian on $C^\infty(X,L^{\otimes N})$ associated with
the connection and the Hermitian structure on $L^{\otimes N}$.
This has the following consequence:
If $\psi_N\in\calH_N\subset L^2(P)$, 
\begin{equation}\label{sobolev2}
\norm{\psi_N}_s = (N^2+N+1)^{s/2}\norm{\psi_N}_0.
\end{equation}
Indeed
elements in $\calH_N$ are eigenfunctions of $\Delta_P$:
We claim that $\Delta_P(\psi_N) = N(N+1)\,\psi_N$.
By virtue of (\ref{thm2ab}), this statement is equivalent to:
$\Delta_h \psi_N = N\psi_N$. This follows from the well-known Bochner-Kodaira
relationship the metric and the $\overline\del$ Laplacian on sections of
$L^N$, see for instance \cite{Duis} Proposition 6.1. 
Clearly (\ref{sobolev2}) has the consequence that (c) implies
(b) for sequences of vectors in $\calH_N$.
\end{proof}

\medskip
We end this subsection with a reminder of the notion of microsupport
in K\"ahler quantization.  Let $\{\psi_N\in\calH_N\}$ be a
sequence of holomorphic sections of the tensor powers of $L$,
and let $\psi = \sum_N\psi_N$.  Since $\Pi(\psi) = \psi$,
the wave-front set of $\psi$ is included in $\calZ$.  The microsupport
of the sequence is defined as the subset of $X$ which is the
projection of $\wf(\psi)$:
we say that $x\in X$ is in the microsupport of $\{\psi_N\}$ if 
and only if
\[
\exists p_x\in P\ \text{such that}\ \pi(p_x) = x\ \text{and}\ 
(p,\alpha_p)\in\wf(\psi).
\]
It follows from the above that microsupport of the sequence is the
empty set iff $\norm{\psi_N}_{L^2} = O(N^{-\infty})$.  In addition,
one can show that $x_0\in X$ is not in the microsupport
of $\{\psi_N\in\calH_N\}$ iff there exists a neighborhood, $V$, of 
$x_0$ such that $\sup_{x\in V}|\psi_N(p_x)| = O(N^{-\infty})$, 
where for all $x$ $p_x\in P$ denotes any point projecting to $x$.

The microsupport has a characterization in terms of the action of
Toeplitz operators analogous to the characterization of the ordinary
wave-front set by the action of pseudodifferential operators.  We refer
to \cite{Ch}, \S 5, for alternative descriptions of the microsupport.

\subsection{Polarized Hermite distributions}
In this section we define and analyze the concept of generalized wave
packets and their symbols in the context of K\"ahler quantization.  
In fact we'll define more general states, associated to isotropic
submanifolds of a quantized K\"ahler manifold (although everything
we do generalizes to any almost K\"ahler manifold quantized by 
a projector, $\Pi$, with the same microlocal structure as the
Szeg\"o projector.)

\bigskip\noindent{\em Definition of polarized Hermite distributions.}

\smallskip

Let us begin by considering a closed conic isotropic submanifold
\[
\calR\subset\calZ.
\]
Obviously $\calR$ is isotropic in $T^*P\setminus\{ 0\}$, and therefore
associated with it are spaces $I^l(P,\calR)$ of Hermite distributions on
$P$.  The general theory of such distributions (together
with many applications) was developed by Boutet de Monvel and
Guillemin, see (see \cite{BG}, or the Appendix for additional remarks). 
The polarized Hermite distributions
associated with $\calR$ are simply the projections of elements of
$I^l(P,\calR)$ by the Szeg\"o projector:

\begin{definition}
The space of polarized Hermite distributions of order $l$
associated with $\calR$ is
\begin{equation}\label{wp0}
I_\Pi^l(P,\calR) := \Pi(I^{l}(P,\calR)).
\end{equation}
\end{definition}
We should point out that the composition Theorem 9.4 of \cite{BG}
one has the inclusion: $I_\Pi^l(P,\calR)\subset I^{l}(P,\calR)$.

\smallskip
Notice that one has a natural isomorphism
\[
\begin{array}{ccc}
P\times\bbR^{+} & \to & \calZ\\
(p,r) & \mapsto & (p, r\alpha_p)
\end{array}
\]
which becomes a symplectomorphism if we put on $P\times\bbR^{+}$
the symplectic structure 
\begin{equation}\label{prsymp}
- d(r\alpha) = -rd\alpha - \alpha \wedge dr
\end{equation}
(here $r$ is the coordinate on the $\bbR^{+}$ factor). 
Using this description of $\calZ$ it is easy to show that the base
of the cone, $\calR$, is a submanifold
$\tilde{Y}\subset P$ such that the infinitesimal
generator of the circle action on $P$ is never tangent to $\tilde{Y}$.
Moreover, the projection, $\pi: P\to X$, restricts to an
isotropic immersion of $\tilde{Y}$ in $X$.  Let us denote by $Y$
this immersed isotropic submanifold of $X$.

\begin{definition}
A sequence of holomorphic sections, $\{\psi_N\in\calH_N\}$, of
the tensor powers of $L$ will be called an {\em Hermite state associated
with $Y$} iff there exists $\psi\in I_\Pi^l(P,\calR)$ such that the
Fourier components of $\psi$ according to (\ref{wp0a}) are precisely
the $\psi_N$.
\end{definition}

For example, if $Y$ is a single point then the coherent states at 
that point are an Hermite state.  The case when $Y$ is Lagrangian
(and hence $\tilde{Y}$ Legendrian) was considered in \cite{BPU}.

\bigskip\noindent{\em Symbolic matters.}

\smallskip

Our next step is to define the symbol of a polarized Hermite
distribution.  We begin by recalling  the nature of the symbol of
general elements in $I^{l}(P,\calR)$.  Such distributions have symbols
which are {\em symplectic spinors} associated with $\calR$.  
The definitions are made in the tangent space to $T^*P$, so for notational 
convenience for each $\rho\in \calR$ we'll set 
$$
R_\rho = T_\rho\calR, \qquad Z_\rho = T_\rho\calZ.
$$
Let
\[
N_\rho := R_\rho^{\circ}/R_\rho
\]
where $R_\rho^{\circ}$ denotes the symplectic orthogonal of $R_\rho$
inside $T_\rho(T^*P)$.  $N_\rho$ is a symplectic vector space, called the
symplectic normal space to $\calR$ at $\rho$.
Abstractly, the symbol of $u\in I^{l}(P,\calR)$ at $\rho$ is a smooth
vector in the metaplectic representation of the metaplectic group
of $N_\rho$ tensored with a half-density along $R_\rho$, i.e.
an element of:
\[
\text{Spin}(R_\rho) := \hffm (R_\rho)\otimes H_\infty (N_\rho).
\]

\begin{lemma}\label{IIn}
Let $E_\rho = \{\,v\in Z_\rho\;;\;\forall u\in R_\rho\ 
\omega(u,v)=0\,\}/R_\rho$ be the symplectic normal of $\calR$ in
$\calZ$ at $\rho$.  Then, the symplectic normal, $N_\rho$, is naturally 
isomorphic to the direct sum
\begin{equation}
N_\rho = E_\rho \oplus Z_\rho^{\circ}
\end{equation}\label{wp1}
where $Z_\rho^{\circ}$ is the symplectic orthogonal of $Z_\rho$ in
$T_\rho(T^*P)$.
\end{lemma}
\begin{proof}
This follows from the fact that
\[
R_\rho^{\circ} = R_\rho^{\circ_\calZ} \oplus Z_\rho^{\circ}
\] 
where $R_\rho^{\circ_\calZ}$ is the symplectic orthogonal of $R_\rho$
inside $Z_\rho$, which itself follows from the fact that $\calZ$ is
a symplectic submanifold of $T^*P$.
\end{proof}

It is important to note that the projection $T_\rho(T^*P) \to T_{\pi(\rho)}X$
induces a symplectic isomorphism $Z_\rho^{\circ} \simeq T_{\pi(\rho)}X$.
On the other hand, because of the negative sign in (\ref{prsymp}),
the same projection takes $E_\rho$ to $(T_{\pi(\rho)}Y^\circ/T_{\pi(\rho)}Y)^-$, 
where the minus indicates a reversal of the symplectic structure.

It follows from Lemma \ref{IIn}
that the metaplectic representation of the metaplectic
group of $N_\rho$ is a tensor product:
\[
H(N_\rho) = H( E_\rho)\hat{\otimes} H(Z_\rho^{\circ})
\]
(Hilbert space tensor product).  It is this decomposition that 
reveals the structure of the symbol of a polarized Hermite distribution.
In order to discuss this structure, we recall that the smooth-vector
factor of the symbol of $\Pi$ is of the form $e\otimes \overline{e}$,
where $e\in H_\infty(Z_\rho^{\circ})$ is a normalized ``ground state"
(which can be identified with the ground state of the harmonic
oscillator on $T_{\pi(\rho)}X$ defined by the metric).

\begin{proposition}\label{WpI}
The symbol of a polarized Hermite distribution, $u\in I_\Pi(P,\calR)$,
is of the form:
\[
\sigma_u = \nu_u\otimes \kappa_u\otimes e,
\]
where $e\in H_\infty (Z_\rho^{\circ})$ is the symbol of the
polarization and
\[
\nu_u\in \hffm (R_\rho)\,\quad \kappa_u \in H_\infty(E_\rho).
\]
By dividing $\sigma_u$ by $e$ one obtains the non-trivial map
in the following exact sequence:
\[
0\to I^{l-1/2}_\Pi(P,\calR)\hookrightarrow I^l_\Pi(P,\calR)\to 
\hffm (R_\rho)\otimes H_\infty(E_\rho)\to 0.
\]
\end{proposition}
We relegate the proof of this technical proposition to an appendix.

\section{Proofs}
\subsection{Proof of Theorem \ref{I}}\label{Pf}

Let $T_f=\{T^\N\,,\,N=1,2,\ldots\}$ be a Berezin-Toeplitz operator
with smooth principal symbol $f:X\to\bbC$.

\bigskip
Part (0) of Theorem \ref{I} is not difficult. 
Without loss of generality we can assume that $\lambda = 0$. 
Since $\{(T^N)^*T^\N$\} is a Toeplitz operator with symbol
$|f|^2$, one has that for all $\psi_N\in\calH_N$,
\[
\norm{T^\N\psi_N}^2 = \inner{(T^N)^*T^\N\psi_N}{\psi_N} \sim 
\inner{\Pi_N(|f|^2\psi_N)}{\psi_N} = 
\inner{|f|^2\psi_N}{\psi_N}.
\]
More precisely, from the definition of B-T operators we have that 
for any sequence of $\psi_N\in\calH_N$,
\[
\norm{T^\N\psi_N}^2 = 
\int_X |f(x)|^2\,|\psi_N(x)|^2\,d\mu_x + \norm{\psi_N}^2\cdot O(1/N),
\] 
where $d\mu$ is the measure on $X$ and $|\psi_N(x)|^2$ is the square of the
length of $\psi_N(x)$ in the Hermitian norm of $L^{\otimes N}\to X$.
It follows that
\[
\inf_{\psi_N\in\calH_N\setminus\{ 0\}}\,
\frac{\norm{(T_f^\N-\lambda I)(\psi_N)}}{\norm{\psi_N}}
\geq \inf_{x\in X}\,|f(x)| + O(1/\sqrt{N}).
\]
In the other direction, let $x_0$ be the point of $X$ where
$\inf |f(x)|$ is attained, and let $\varphi_{p_0}^{(N)} = \Pi(\cdot,p_0)$
be a coherent
state at a point $p_0\in P$ that projects down to $x_0$.  Then
\[
\inf_{\psi_N\in\calH_N\setminus\{ 0\}}
\frac{\norm{T^\N\psi_N}^2}{\norm{\psi_N}^2} \leq
\frac{\inner{(T^N)^*T^\N\varphi_{p_0}}{\varphi_{p_0}}}{\norm{\varphi_{p_0}}^2}
= |f(p_0)|^2 + O(1/N).
\]

\medskip
We now prove part (1) of the Theorem.  
We begin with a few preliminary considerations.
By Proposition 2.13 of \cite{BG}, there exists a classical
pseudodifferential operator of order zero on the circle bundle $P$,
$Q$, such that:
\begin{enumerate}
\item $Q$ commutes with the Szeg\"o projector ($[\Pi\,,\,Q] = 0$) and with
the $S^1$ action.
\item For each $N$ the restriction of $Q$ to $\calH_N$,
$ Q: \calH_N\to\calH_N$ is equal to $T^\N$.
\item The principal symbol of $Q$ satisfies:
\[
\forall (p,r\alpha_p)\in\Sigma\qquad \sigma_Q(p,r\alpha_p) = f(\pi(p))
\]
where $\pi: P\to X$ is the projection.
\end{enumerate}

Assume now that $\lambda = f(x_0)$, where $x_0\in X$ is such that
\begin{equation}\label{3a}
\PB{\Re f}{\Im f}(x_0) < 0.
\end{equation}
The inverse image $f^{-1}(\lambda)$ is, in a neighborhood of $x_0$, a
codimension-two symplectic submanifold of $X$.  Let us pick an
isotropic submanifold of $f^{-1}(\lambda)$ (not necessarily closed),
$Y$, containing $x_0$.  Then $Y$ is an isotropic submanifold of $X$.
Our considerations are local:  we restrict our attention to a
neighborhood of $x_0$ where $\PB{\Re f}{\Im f}$ is negative and such
that there exists a lift of $Y$ to a conic isotropic submanifold,
$\calR\subset\calZ$, in the sense of the previous section.
We will construct a pseudomode with microsupport equal to $Y$.
Notice that we may take $Y=\{x_0\}$ if we wish.

Since $[\Pi, Q] = 0$ and $\Pi$ is self-adjoint,
we have $[\Pi, Q^\dagger] = 0$.  Therefore, by
Proposition 11.4 of \cite{BG},
\begin{equation}\label{3b}
\forall (p,r\alpha_p)\in \Sigma
\qquad \PB{\Re \sigma_Q}{\Im \sigma_Q} (p,r\alpha_p) =
\PB{\Re f}{\Im f} (\pi (p)),
\end{equation}
where $\sigma_Q$ is the principal symbol of $Q$.  Notice that the
Poisson bracket on the left is on $T^*P$ (with respect to the
cotangent bundle structure), while the one on the right is the
Poisson bracket on $X$.  Therefore, the Poisson bracket conditions
(\ref{3a}, \ref{3g}) on $f$ are inherited by $\sigma_Q$.

To see how the Poisson bracket condition becomes relevant in our
considerations, we first consider a calculation in the Heisenberg
representation on $L^2(\bbR^k)$.
\begin{lemma}\label{heislemma}
Let $\calL$ be the operator on $L^2(\bbR^k)$ corresponding to the action
of $\xi \in (\bbR^{2k},\omega) \otimes \bbC$ under the Heisenberg
representation.  Then if
$$
\omega(\Re \xi, \Im\xi) > 0,
$$
then the restriction of $\calL$ to the smooth vectors maps $\calS(\bbR^k)$
onto itself, with a non-zero kernel.
\end{lemma}
\begin{proof}
The metaplectic representation describes how $\calL$ transforms
under the action of the symplectic group on $\xi$.  That is, $\calL_{g.\xi}
= U(g) \calL_\xi U(g^{-1})$, where $g\mapsto U(g)$ is the projective unitary
representation that gives rise to the metaplectic representation when we
take the double cover.  The action of the metaplectic group preserves the
smooth vectors $\calS(\bbR^k)$, so in our argument we can replace $\xi$ by
$g.\xi$ for $g$ symplectic.

Under the assumption that $\epsilon := \omega(\Re \xi, \Im\xi) > 0$, it is a
straightforward  exercise to see that $g$ can be chosen
so that $g.\xi = \epsilon e_1 + i f_1$, where $e_1,\dots,e_k,
f_1,\dots,f_k$ is the standard symplectic basis for
$\bbR^{2k}$.
Thus it suffices to prove the result for
$$
\calL = \frac{\partial}{\partial x_1} + \epsilon x_1.
$$

We see then that $\ker\calL \cap \calS(\bbR^k)$ contains functions of the
form  $\psi(x_1,\dots,x_k) = e^{-\epsilon x_1^2/2}  a(x_2,\dots, x_k)$.
And to show that $\calL$ maps $\calS(\bbR^k)$ onto itself, let $f\in
\calS(\bbR^k)$.  The ODE $\calL u = f$ can be solved by variation of
parameters:
$$
u(x_1,\dots,x_k) = \int_0^{x_1} f(t,x_2,\dots,x_k) e^{\epsilon(t^2-x_1^2)/2}
\>dt.
$$
For $\epsilon>0$, the estimates showing that $u\in\calS(\bbR^k)$ follow
easily.
\end{proof}

Applied at the symbol level, Lemma \ref{heislemma} leads directly to the
following construction:
\begin{proposition}
There exists a distribution $u$ in the class $I_\Pi^{0}(P,\calR)$ of
polarized Hermite distributions of order zero associated with $\calR$
such that:
\begin{equation}\label{3c}
(Q-\lambda I)u \in C^\infty(P).
\end{equation}
\end{proposition}
\begin{proof}
Suppose $u\in I_\Pi^{0}(P,\calR)$ with $\sigma_u = \nu_u \otimes \kappa_u \otimes e$
as in Proposition \ref{WpI}.  It's clear that
$(Q-\lambda I)u = \Pi(Q-\lambda I)u$  is also a polarized Hermite distribution.  What is
its symbol? Since the principal symbol of $Q-\lambda I$ vanishes on $\calR$,
we are led to use the first transport equation for the Hermite calculus.
The Hamilton vector field, $\xi$, of the
symbol of $Q$ at a point $\rho\in\calR$ is in the symplectic normal
space $N_\rho$.  Therefore, the Heisenberg representation of that
space associates to $\xi_\rho$ an operator, $\calL$, on the space
$H_\infty(N_\rho)$.  According to Theorem 10.2 in \cite{BG},
$(Q-\lambda I)(u)\in I^{-1/2}(P,\calR)$ and its symbol
is $\nu_u \otimes \calL(\kappa_u\otimes e)$.  
In fact, since $(Q-\lambda I)u$ is still polarized,
under the decomposition $H_\infty(N_\rho) \simeq
H_\infty(E_\rho) \otimes H_\infty(Z^\circ)$, $\calL$ acts only on 
$H_\infty(E_\rho)$.  The symbol is really $\nu_u\otimes \calL(\kappa_u) \otimes e$.

We noted in the previous section that the pull-back of the symplectic form from $X$ 
under the natural projection $E_\rho \to T_{\pi(\rho)}X$ is the opposite of the 
symplectic form on $E$. 
So the Poisson bracket condition (\ref{3b}) along with Lemma \ref{heislemma}
implies that, as an operator on $H_\infty(E_\rho)$, $\calL$ is onto and
has a non-trivial kernel.

Because $\calL$ has a kernel, we can choose 
$u_0\in I^{0}(P,\calR)$ with symbol $\calL(\sigma(u_0)) = 0$.    
Therefore $v_1:=(Q-\lambda I)(u_1)\in
I^{-1}(P,\calR)$.   Because $\calL$ maps onto $H_\infty(E_\rho)$ 
(by Lemma \ref{heislemma} again), we can then
find $u_1\in I^{-1/2}(P,\calR)$ such that
$\calL(\sigma_{u_1}) = - \sigma_{v_1}$, 
Thus $v_2 :=(Q-\lambda I)(u_0+u_1)\in I^{-3/2}(P,\calR)$.
Continuing in this fashion and finishing with a Borel summation of the
$u_j$'s produces $u \in I^{0}(P,\calR)$ such that $(Q-\lambda I)(u)$ is of
order $(-\infty)$, and therefore smooth.
\end{proof}

\bigskip
We can now finish the proof of the part (1) of Theorem \ref{I}.  
Let $u$ be as in the previous Proposition, which we choose to
have a Gaussian principal symbol.
Define $u_N:=\Pi_N(u)$.
By Lemma \ref{Sobolev}, (\ref{3c}) implies the norm estimates
\begin{equation}\label{3e}
\norm{(T^\N-\lambda I)u_N} = O(N^{-\infty}).
\end{equation}
On the other hand, letting $\varphi_{p_0}^{(N)}$ be
the coherent state at $p_0\in P$, from the reproducing property $u_N(p_0)
= \langle \varphi_{p_0}^{(N)}, u_N \rangle$ we obtain the estimate
$$
\norm{u_N} \ge \frac{|u_N(p_0)|}{\norm{\varphi_{p_0}^{(N)}}}.
$$
Combining (\ref{psirho}) with the asymptotics 
$$
\norm{\varphi_{p_0}^{(N)}} = \Pi_N(p_0,p_0) \sim \Bigl(\frac{N}{2\pi}\Bigr)^n
$$
(see e.\ g.\ equation (28) in \cite{Zel}),
we see that 
$$
\norm{u_N} \ge C N^{-(l+1)/2},
$$
where $l = \dim Y$.  Together with (\ref{3e}) this implies (\ref{1g}).

\bigskip
To prove part (2) of Theorem \ref{I}, we start by assuming
that
\begin{equation}\label{3g}
\PB{\Re f}{\Im f}(x_0) > 0,
\end{equation}
and let $u_N\in\calH_N$ be a sequence of vectors such that 
\[
\norm{(T_f^\N-\lambda I)(u_N)} = O(N^{-\infty})
\quad\text{and}\quad \norm{u_N} = 1.
\]
If we let $u$ be the distribution on $P$ whose Fourier coefficients are the
$u_N$, then we can rewrite this condition as
\begin{equation}\label{3h}
(Q-\lambda I)(u) \in C^\infty(P).
\end{equation}
We will now quote Theorem 27.1.11. of \cite{Ho4}  
(H\"ormander, Vol IV) asserting that the Poisson
bracket condition (\ref{3g}), translated to the corresponding statement about 
$\sigma_Q$, implies that $Q-\lambda I$ is microlocally subelliptic on the set 
\[
\Sigma_{x_0} := \{ \rho\in\calZ ;\; \pi(\rho) = x_0\}\subset T^*P,
\]
with loss of $1/2$ derivatives.  Therefore, by (\ref{3h}),
for all $s\in\bbR$, $u\in H_{(s)}^{\text{\tiny loc}}$
at every $\rho\in\Sigma_{x_0}$.  What this means is that
for each such $\rho$ we can write: $u = u_0+u_1$
where $u_1\in H_{(s)}(P)$ and $\rho\not\in\text{WF}(u_0)$,
for all $s$.  Therefore the wave-front set of $u$ is disjoint
from $\Sigma_{x_0}$.

\subsection{Proof of Theorem \ref{II}}
There are two ingredients in the proof, one is a very general localization
statement (well-known in the theory of $\h$-admissible $\Psi$DOs),
and the second H\"ormander's results on microlocal subellipticity
that we used in the previous section.  We begin with the
localization result, which is of interest in its own right. 

\begin{proposition}\label{IIb}
Let $T_f^\N$ be a Berezin-Toeplitz operator with principal symbol $f$,
and
let $\{u_N\in\calH_N\}$ be a minimizing sequence of the Rayleigh
quotients,
$\{\frac{\norm{(T_f^\N-\lambda I)(\psi)}}{\norm{\psi}},\psi\in\calH_N\}$, 
where $\norm{u_N}=1$ for all $N$.  Then the distribution
$u:=\sum_N u_N\in C^{-\infty}(P)$
has wave-front set contained in the set of
points, $\rho\in T^*P$, such that $\pi(\rho)\in f^{-1}(\lambda)$, and
therefore the microsupport of $\{u_N\}$ is contained in $f^{-1}(\lambda)$.
\end{proposition}
\begin{proof}
For simplicity of notation, assume without loss of generality
that $\lambda = 0$.  Let $Q$ denote the operator on $P$
inducing the $T^\N_f$ and commuting with $\Pi$, as before.
The minimizing sequence is a sequence of
eigenstates of the non-negative, self-adjoint 
classical pseudodifferential operator of order zero
$S = Q^*Q$ on $P$.  We will denote by $S^\N:\calH_N\to\calH_N$ the 
restriction of $S$ to $\calH_N$.  By assumption $u_N$ is an
eigenvector of $S^\N$ corresponding to the smallest eigenvalue.
Let $\chi\in C_0^\infty(\bbR)$ be a test function, $R$ 
a zeroth-order $\Psi$DO on $P$, and consider the traces
\[
\Upsilon_{\chi,R}^\N = \tr \int\ \chi(t)\,R\,e^{-itNS_N}\,dt.
\]
If we write the eigenvalues and eigenvectors of $S_N$ in the form:
\[
S_N(\psi_j^\N) = E_j^\N\,\psi_j^N,\quad E_1^\N\leq E_2^\N\leq\cdots\
E_{d_N}^\N
\]
where $d_N = \dim \calH_N$ and $\{\psi_j^N\}_j$ is an orthonormal
basis of $\calH_N$, then
\[
\Upsilon_{\chi,R}^\N = \sum_{j=1}^{d_N}\,\hat\chi(NE_j^\N)\,
\inner{R\psi_j^\N}{\psi_j^\N}.
\]
Just as in the proof of part (0) of Theorem \ref{I},
by taking coherent states as trial functions one
obtains the estimate: $E_1^\N = O(1/N)$.

\begin{lemma} If the microsupport of $R$ is
disjoint from the characteristic set of $Q$, then
$\Upsilon_{\chi,R}^\N = O(N^{-\infty})$.
\end{lemma}
\noindent{\em Proof of the Lemma:}  
The operator $e^{-itD_\theta S}$ is a Fourier integral operator,
and a simple wave-front set calculation shows that the
wave-front set of the operator
\[
\Pi\circ S_\chi := \Pi\circ \int \chi(t)\,e^{-itD_\theta S}\,dt
\]
is contained in the set
\[
\{(\rho,\rho')\in\calZ\times\calZ\;;\; \sigma_Q(\rho) = 0 = \sigma_Q(\rho')\}.
\]
Therefore, if the microsupport of $R$ is disjoint from the characteristic
set of $Q$, the operator $R\circ\Pi\circ S_\chi$ is smoothing.
Consider next the generating function of the $\Upsilon^\N$,
\[
\Upsilon (s) := \sum_N e^{iNs}\,\Upsilon_{\chi,R}^\N =
\Tr \Bigl(U(e^{is})\circ R\circ\Pi\circ S_\chi\Bigr),
\]
where $U(e^{is})$ is the operator on $P$ given by composition 
by the action of $e^{-is}$.
Another wave-front set calculation shows that
$\Upsilon\in C^\infty(S^1)$ because $R\circ\Pi\circ S_\chi$ is
smoothing, and therefore the Fourier coefficients of $\Upsilon$ 
are rapidly decreasing.  This proves the lemma

\medskip
Continuing with the proof of the Proposition, choose $\chi$ so that
$\hat\chi > 0$ (we can even take it so that $\hat\chi$ is equal to
one in a neighborhood of zero),  and choose $R$ of the form:
$ R = F^*\circ F $ where $F$ is a zeroth order $\Psi$DO on $P$. 
Then all terms in the sum defining $\Upsilon_{\chi,R}^\N$ are
non-negative, and the previous Lemma implies that, if the 
microsupport of $F$ is disjoint from Char$(Q)$, one has:
$ \norm{F(u_N)} = O(N^{-\infty})$.  But this implies that, for any 
such $F$,
\[
F(u) = \sum_N F(u_N)\in C^\infty (P).
\]
Since we can pick $F$ microlocally elliptic in a neighborhood of
any point in the complement of Char$(Q)$,
the wave-front set of $u$ must be contained in Char$(Q)$.
\end{proof}

Turning to the proof of Theorem \ref{II},
we note that under its assumptions we are in a position to apply
the subellipticity results of
H\"ormander.  Note first that any given
repeated Poisson bracket of the real and imaginary parts of $f$
evaluated at a point $x\in X$ is equal to the same repeated 
Poisson bracket of the real and imaginary parts of the
principal symbol of $Q$, evaluated at any $\rho\in\calZ$ such that
$\pi(\rho) = x$.  


Next, we claim that the hypotheses of Theorem
\ref{II} imply the hypotheses of Theorem 27.1.11. of \cite{Ho4},
namely:
\begin{itemize}
\item[(A)] For every $x\in f^{-1}(\lambda)$ there is some $j\leq k$ and
some $z\in\bbC$,
\begin{equation}\label{27}
\Bigl((\Xi_{\Re (zf)})^j\,\Im (zf)\Bigr) (x) \not=0
\end{equation}
where $\Xi_{\Re (zf)}$ is the Hamilton vector field of $\Re (zf)$,
considered as a differential operator.
\item[(B)] The repeated Poisson bracket
above is non-negative (it is positive, actually)
if $j$ is the smallest integer such that (\ref{27}) holds for
some $z$.  Moreover such $j$ is odd.
\end{itemize}
Part (A) follows from hypothesis (2) in Theorem \ref{II} by
Corollary 27.2.4 in \cite{Ho4}.  Lemma 5.1
of \cite{DSZ} shows that hypothesis (1) of Theorem \ref{II} implies
H\"ormander's Condition ($\overline\Psi$), which, as indicated in the 
remark
after Theorem 27.1.11 of \cite{Ho4}, implies (B).  The proof of this 
remark
is further detailed in the first paragraph in the proof of Theorem 4 of 
\cite{DSZ}.

\medskip
Let $\Theta\subset\calZ$ be the set of points in $\calZ$
projecting to $f^{-1}(0)$ (recall that we took $\lambda = 0$).
We conclude, by Theorem 27.1.11 of \cite{Ho4} and the previous
considerations, that for all $\rho\in\Theta$,
$Q$ is subelliptic at $\rho$ with loss of at most $\delta = k/(k+1)$
derivatives.

\medskip
By Lemma 27.1.5 of \cite{Ho4}, for each $\rho\in \Theta$
there exists a zeroth-order $\Psi$DO, $A_\rho$, non-characteristic
at $\rho$ and such that
\[
\forall g\in C^\infty(P)\quad
\norm{A_\rho g}_{(-\delta)} \leq C_\rho
\Bigl( \norm{Qg}_{L^2} + \norm{g}_{(-1)} \Bigr).
\]
An examination of the proof of this lemma shows that we can 
take the symbol of each $A_\rho$ to be non-negative.
Since $\Theta$ is a cone with compact base, there exists
an integer $K$ such that the sum of $K$ of the operators
$A_\rho$, call them $A_1,\ldots, A_K$, is non-characteristic
at each $\rho\in\Theta$.  We denote such a sum
by $A = \sum_{j=1}^K A_j$.  Then
\begin{equation}\label{thm5}
\forall g\in C^\infty(P)\quad
\norm{Ag}_{(-\delta)} \leq C
\Bigl( \norm{Qg}_{L^2} + \norm{g}_{(-1)} \Bigr)
\end{equation}
for some fixed $C>0$.  By averaging with respect to the $S^1$
action on $P$ (which preserves Sobolev norms), we can further
assume without loss of generality that $[A,D_\theta]=0$.

Let $B$ be a microlocal parametrix of $A$ in a neighborhood of $\Theta$
such that $[B,D_\theta]=0$.
Since, by Proposition \ref{IIb}, the wave-front set of
$u$ is contained in $\Theta$, we have:
\[
u = BA(u)+g,\quad \text{with}\ g\in C^\infty(P).
\]
By taking Fourier components it follows that
\[
\forall N\quad \norm{u_N}_{(-\delta)} \leq C_B
\norm{A(u_N)}_{(-\delta)} + \norm{g_N}_{(-\delta)},
\]
where $C_B$ is the $H_{(-\delta)}$ norm of $B$.
This together with (\ref{thm5}) implies that for all $N$
\[
\norm{Q(u_N)}_{L^2} \geq C_1 \norm{u_N}_{(-\delta)} -
C_2 \norm{g_N}_{(-\delta)} - \norm{u_N}_{(-1)}.
\]
Since $g\in C^\infty(P)$, $\norm{g_N}_{(-\delta)} = O(N^{-\infty})$,
while equation (\ref{sobolev2}) gives us that
\[
C_1 \norm{u_N}_{(-\delta)} -\norm{u_N}_{(-1)} = 
C_1(N^2+N+1)^{-\delta/2} - (N^2+N+1)^{-1/2}
\]
Since $\delta = k/(1+k)<1$, this proves Theorem \ref{II}.

\section{Examples}
\newcommand{\SU}{\text{SU}}
\newcommand{\CP}{{\bbP ^1}}

We now look at specific examples of quantized K\"ahler manifolds,
$X$, and of Toeplitz operators.  The corresponding Hilbert spaces have
canonical bases and therefore the Toeplitz operators become sequences
of matrices of a specific type that we compute.

\subsection{A preliminary example}

We begin with a concrete example associated with $X=\bbC$, the plane
with its usual complex structure.  Although this $X$ is not compact
(and therefore it does not fit precisely into the general framework of
this paper) we will ``cut it" to the unit disk, both symplectically and
quantum-mechanically.  This leads to a sequence of large matrices to
which microlocal methods apply, provided one stays away from the
boundary of the unit disk.  We only consider an explicit operator which is
the microlocal model of the general case.  It will be clear that what
we do easily generalizes to other operators in this setting.
In this section we want to be explicit and avoid using the general machinery.

Recall that the K\"ahler quantization of the plane gives rise to the
Bargmann spaces
\[
\calB_N = \{\, f:\bbC\to\bbC\ \text{entire}\;;\;
\norm{f}^2:=\frac{1}{\pi}\,\int_X |f(z)|^2\,e^{-Nz\zbar}\,dxdy <\infty \}
\]
where $z=x+iy$ and $N>0$.  Elements of $\calB_N$ arise
from the general K\"ahler quantization scheme applied to $\bbC$.
The quantizing line bundle $L\to\bbC$ is holomorphically trivial,
and so its sections can be identified with entire functions on $\bbC$.
The Hermitian structure on $L$, however, is not trivial.
We introduce the following notation for the length function
of $\psi\in\calB_N$ as a section of $L^N\to\bbC$:
\begin{equation}\label{prelim}
|\psi(z)|_s := |\psi(z)|\,e^{-Nz\zbar/2}.
\end{equation}
Notice that then the norm of $\psi$ is the integral of the
function $|\psi|_s$ with respect to the area form.

A fundamental operator on $\calB_N$
is the harmonic oscillator (shifted by 1/2, for convenience), 
\[
\text{Op}(H_{\text{\tiny H.O.}}) = N^{-1} z\frac{d\ }{dz}
\]
which is a Berezin-Toeplitz operator with symbol
$H_{\text{\tiny H.O.}}(z,\zbar) = z\zbar$.  The eigenfunctions and
eigenvalues of $\text{Op}(H_{\text{\tiny H.O.}})$ are:
\[
N^{-1} z\frac{d\ }{dz} z^j = N^{-1} j z^j,\qquad \text{and}\quad
\norm{z^j} = N^{-1(j+1)/2}\sqrt{j!},
\]
so that
\[
\ket{k} = \frac{N^{(k+1)/2}}{\sqrt{k!}}\,z^k,\qquad k=0,1,\ldots 
\]
is an orthonormal basis of $\calB_N$.
The unit disk is the region of phase space where the classical
energy $ z\zbar$ is less than one.  The analogous object
quantum-mechanically is the span of the eigenfunctions with eigenvalue
less than one, that is the monomials $z^j$ with $j\leq N$.  Thus
\[
\calH_N = \{\,\text{polynomials in the complex variable }z
\text{ of degree }\leq N\;\}
\]
where it is now natural to restrict $N$ to be an integer. 
This setting is close to the case of the
sphere to be considered in the next section, provided one does not get
too close to the boundary of the unit disk.  The
Hilbert space $\calH_N$ of the sphere can also be identified with the
space of polynomials in a complex variable of degree at most $N$, although
on the sphere $\norm{z_j}$ is essentially $(C^j_N)^{-1/2}$.
Symplectically, the sphere is the disk with its boundary collapsed to a point.

\renewcommand{\P}{\Theta^{\bot}}
The example we will be studying is based in the following observation.
Suppose one has a Berezin-Toeplitz
operator, $Q$, and a state $\psi$ such that $Q\psi=0$. 
Suppose one has a ``good" semi-classical
cut-off operator, $\Theta_N$ (a projector), which is semi-classically the identity
in a certain region of the plane (referred to as the allowed region),
and let $P_N = \Theta_N Q\Theta_N$.  Then
\[
0 = \Theta_N Q\Theta_N(\psi) + \Theta_N Q\P_N(\psi) + \P_N Q\Theta_N(\psi) +
\P_N Q\P_N(\psi)
\]
where $\P_N = I-\Theta_N$ is the complementary projection.  
The second and third terms in this sum will be very small due to the
assumed localization properties of $\Theta_N$.
It follows that $\Theta_N(\psi)$ is a good pseudo-mode
for $P_N$ (with pseudo-eigenvalue zero) if
\[
\frac{\norm{\P_N\psi}}{\norm{\Theta_N\psi}}
\]
is small, that is, if $\psi$ concentrates in the classically allowed region
of $\Theta_N$.  We will take $\Theta_N$ to be the orthogonal projection 
\[
\Theta_N:\ \calB_{N}\to \calH_N,
\]
for which the classically allowed region is
the interior of the unit disk.  We now proceed to make this statement 
more precise.  Our basic tool is the reproducing kernel of $\calB_N$,
in the form of the coherent states: For each $w\in\bbC$, let 
\begin{equation}\label{prelim0}
\varphi_w(z) := N e^{Nz\overline{w}}\in\calB_N.
\end{equation}
These states have the reproducing property
\begin{equation}\label{prelim0a}
\forall\psi\in\calB_N,\ z\in\bbC\qquad \psi(z) = \inner{\psi}{\varphi_z}.
\end{equation}
Notice that, in particular
\begin{equation}\label{prelim0b}
\norm{\varphi_w}^2 = \inner{\varphi_w}{\varphi_w} = N\,e^{N|w|^2}.
\end{equation}

\begin{lemma}  For every $\delta>0$ there exists $n>0$
such that for all $N>n$ and for each $w\in\bbC$ such that $|w|<1$,
we have:
\begin{equation}\label{prelim1}
\norm{\Theta_N^\bot \varphi_w}^2 \leq \frac{1+\delta}{\sqrt{2\pi}}
\norm{\varphi_w}^2\,|w|^2\,
\sqrt{N}\,e^{-N(1-|w|^2)^2/2}.
\end{equation}
\end{lemma}
\begin{proof}
From the Taylor series expression
for $\varphi_w$ and the orthogonality relations of the monomials $z^j$, one
can show that 
\begin{equation}\label{prelima}
\norm{\Theta_N^\bot \varphi_w}^2 = \norm{\varphi_w}^2\,
\Bigl( 1-\frac{\Gamma(N+1,N|w|^2)}{N!}\Bigr),
\end{equation}
where $\Gamma(n,x) = \int_x^\infty t^{n-1}e^{-t}dt$
is the incomplete gamma function.
To estimate the quantity in parenthesis, notice that
\[
N! = \Gamma(N+1) = \Gamma(N+1,N|w|^2) + \int_0^{N|w|^2}\,e^{-t}\,t^N\,dt.
\]
Therefore, dividing by $N!$ and making the change of
variables: $s=t/N$, we see that the quantity in parenthesis in
(\ref{prelima}) equals
\[
\frac{N^{N+1}}{N!}\,\int_0^{|w|^2}\,e^{-Ns}\,s^N\,ds = 
\frac{N^{N+1}}{N!}\,e^{-N}\,\int_0^{|w|^2}\,e^{-Nf(s)}\,ds
\]
where $f(s) = s-\log(s)-1$.  This is a decreasing function on 
$s\in(0,1)$, and therefore the last integrand is maximal 
at $s=|w|^2$.  On the other hand, it is elementary to check that
\[
\forall\epsilon\in (0,1)\quad f(1-\epsilon) \geq \epsilon^2/2,
\]
and therefore
\[
\frac{N^{N+1}}{N!}\,\int_0^{|w|^2}\,e^{-Ns}\,s^N\,ds \leq
\frac{N^{N+1}}{N!}\,e^{-N}\,|w|^2\,e^{-N(1-|w|^2)^2/2}.
\]
An application of Stirling's formula finishes the proof.
\end{proof}

As we now see this Lemma implies the localization properties
of $\Theta_N$:
\begin{corollary}
For each $\psi\in\calB_{N}$, and for each $z$ such that $|z|<1$,
\begin{equation}\label{prelim3}
|\Theta_N^\bot(\psi)(z)|_s \leq   \norm{\psi}\,|z|\,
N^{3/4}\,e^{-N (1-|z|^2)^2/4}.
\end{equation}
Therefore, for all $\epsilon >0$ there exist $C,\ a>0$ such that
for all $\psi\in\calB_{N}$
\begin{equation}\label{prelim3.0}
\dint_{|z|<1-\epsilon} |\Theta_N^\bot(\psi)(z)|^2_s\,dx\,dy \leq
C\,\norm{\psi}^2\,N^{3/2}\,e^{-aN}.
\end{equation}
\end{corollary}
\begin{proof}
By the reproducing property (\ref{prelim0a}), we have:
\[
\Theta_N^\bot(\psi)(z) = \inner{\Theta_N^\bot\psi}{\psi_z}=
\inner{\psi}{\Theta_N^\bot\varphi_z}, 
\]
and therefore, by the Cauchy-Swartz inequality and (\ref{prelim1}),
for all sufficiently large $N$
\[
|\Theta_N^\bot(\psi)(z)| \leq   \norm{\psi} \norm{\varphi_z}\,|z|\,
N^{1/4}\,e^{-N(1-|z|^2)^2/4}.
\]
Using (\ref{prelim0b}) we obtain (\ref{prelim3}), and (\ref{prelim3.0})
follows by integration.
\end{proof}

We will need one more general fact about the projection, $\Theta_N$:

\begin{lemma}
For all $z\in\bbC$ and $\psi\in\calH_{N}$, one has:
\begin{equation}\label{prelim3a}
|(\Theta_N\psi)(z)|_s \leq (N+1)\,\max_{0\leq t\leq 2\pi}\,|\psi(e^{it}z)|_s.
\end{equation}
\end{lemma}
\begin{proof}
This is a consequence of the following formula for $\Theta_N$:
\begin{equation}\label{prelim3b}
\Bigl(\Theta_N\psi\Bigr)\, (z) =\frac{1}{2\pi}\, \sum_{k=0}^N\,\int_0^{2\pi}\,
e^{-ikt}\,\psi(e^{it}z)\,dt,
\end{equation}
together with the fact that the Hermitian weight function $e^{-Nz\zbar/2}$
is $S^1$ invariant.
\end{proof}

\newcommand{\phiz}{\varphi_{\mu,z_0}}
\bigskip
We now show how the above results can be applied to to a concrete example,
the sequence $P_N  = \Theta_N Q\Theta_N$ where $Q$ is the model operator: 
\[
Q = \frac{1}{N}\frac{d}{dz} + \mu z 
\]
acting on the space $\calH_N$.  The symbol of $\{P_N\}$ is the function
$f=\mu z+ \zbar$ restricted to the unit disk.  Notice that
$\PB{\Re f}{\Im f}$ is identically equal to $\mu^2-1$, which we will
assume is negative, i.e.\ we now take $\mu\in (-1,1)$.  
Notice that for each $\lambda\in\bbC$ there exists a unique $z_0$ such that 
$\lambda = f(z_0)$.

Every complex number is an eigenvalue of the operator $Q$.  Specifically,
for each $z_0\in\bbC$ consider the state
\[
\phiz (z) := \sqrt{N}\,e^{-N |z_0|^2/2}\,e^{Nz\overline{z_0}}\,
e^{-\frac{N}{2}\mu (z-z_0)^2}.
\]
One can verify that $\phiz\in\calB_{N}$ because $|\mu|<1$.
The state $\phiz$ is the (quantum) translate of
the basic ``squeezed state" at the origin, $\sqrt{N}e^{-\frac{N}{2}\mu z^2}$,
to the point $z_0$, and its norm is a universal constant (independent of
$N$ and of $z_0$).  It is trivial to verify that
\[
Q\phiz(z) = f(z_0)\,\phiz(z).
\]
As a section of $L^{\otimes N}\to\bbC$, the length of $\phiz$ at $z$ is
\begin{equation}\label{prelimz}
|\phiz(z)|_s =  \sqrt{N}\,e^{-N F(z)}
\end{equation}
where $F$ is the real quadratic function
\[
F(z) = -\Re\Bigl(z\overline{z_0} -\frac{1}{2}\mu (z-z_0)^2 \Bigr) +
z\zbar/2 + |z_0|^2/2.
\]
$F$ is non-negative,
vanishing exactly at $z=z_0$.  The states $\phiz$ are the
model of the polarized Hermite states associated with a point.
Our result is as follows:

\begin{proposition}
With the previous notation and if $\mu\in(-1,1)$ and $|z_0|<1$,
\[
\frac{\norm{(P_N-\lambda)\Theta_N(\phiz)}}{\norm{\Theta_N(\phiz)}} =
O(N^{1/4}\,e^{-aN})
\]
for some $a>0$.
\end{proposition}
\begin{proof}
We first will show that
\begin{equation}\label{prelim4z}
\norm{\Theta_N(\phiz)}^2= \norm{\phiz}^2 + O(e^{-aN}). 
\end{equation}
As noted above, $\norm{\phiz}$ is a constant independent of $N$. 
For simplicity, $a$ will
denote a positive constant that may not be the same at each occurrence.

Let $\Delta$ denote a disk
of radius less than one containing $z_0$ in its interior.
Then, by (\ref{prelim3a}) and the decay properties of $\phiz$
\begin{equation}\label{prelim4a}
\norm{\Theta_N(\phiz)}^2 = \dint_{\Delta} |\Theta_N(\phiz)|^2_s\,dx\,dy +
O(e^{-aN}).
\end{equation}
Next, notice that on  $\Delta$,
$|\phiz(z)|_s$ is uniformly bounded by a constant times $\sqrt{N}$,
and therefore (by (\ref{prelim3a}) again)
$|\Theta_N(\phiz)(z)|_s$ is bounded by a constant times
$N^{3/2}$ there.  It follows that
\[
\Bigl|\, |\Theta_N(\phiz)(z)|^2_s - |\phiz(z)|_s^2\,\Bigr| = 
|\Theta_N^\bot(\phiz)(z)|_s\cdot O(N^{3/2}),
\]
with a constant uniformly on $z\in\Delta$.  Therefore, by
(\ref{prelim3.0}),
\[
\dint_{\Delta} |\Theta_N(\phiz)|^2_s\,dx\,dy =
\dint_{\Delta} |\phiz|^2_s\,dx\,dy + O(e^{-aN}).
\]
Once again, by the localization properties of $\phiz$ we have
\[
\norm{\phiz}^2 = \dint_{\Delta} |\phiz|^2_s\,dx\,dy +
O(e^{-aN}).
\]
The last two equations and (\ref{prelim4a})
imply (\ref{prelim4z}).

Let us now turn our attention to the vector
$(P_N-\lambda)\Theta_N(\phiz)$.  In the standard orthonormal
basis $\{\ket{k}\}$, the matrix $T=(t_{lm})$ of $P_N$ is tri-diagonal.
Specifically, the only non-zero elements of this matrix are:
\begin{equation}\label{prelimm}
t_{k,k+1} = \sqrt{k/N} \quad\text{and}\quad t_{k+1,k} = \mu\sqrt{k/N}.
\end{equation}
Let $\phiz  = \sum_{k=0}^\infty a_k\ket{k}$.
Since $\phiz$ is an eigenfunction of $Q$ with eigenvalue $\lambda$, 
\[
(P_N-\lambda)\Theta_N(\phiz) = \Bigl(\mu\, a_{N-1} -\lambda\,a_N\Bigr)\ket{N}
\]
We claim that both $a_{N-1}$ and $a_N$ are exponentially small in $N$.
Indeed
\[
a_N\ket{N} = \frac{1}{2\pi}\,\int_0^{2\pi}\,e^{-iNt}\,\phiz(e^{it}z)\,dt,
\]
which implies that for all $z\in\bbC$
\[
|a_N\,z^N|\leq \frac{\sqrt{N!}}{N^{(N+1)/2}}\,\max_t |\phiz(ze^{it})|,
\]
where the absolute value is the standard one.  Evaluating both
sides at $z=1$ and applying Stirling's formula we obtain
\[
|a_N| \leq C N^{-1/4}\,e^{-N/2}\,
\max_t |\phiz(e^{it})|.
\]
But $ e^{-N/2}\max_t |\phiz(e^{it})| = \max_t |\phiz(e^{it})|_s$,
and therefore
\[
e^{-N/2}\max_t |\phiz(e^{it})| = O(\sqrt{N}e^{-aN})
\]
since $|z_0|<1$ (see \ref{prelimz} and the remarks following it). 
Therefore $a_N = O(N^{1/4}e^{-aN})$, and similarly for $a_{N-1}$.
\end{proof}

Much more generally we can start with a pseudo-mode of a general
Berezin-Toeplitz operator, $Q$, on Bargmann space. By the localization
properties of the projector $\Theta_N$, the projection by $\Theta_N$ of the
pseudomode will be a pseudomode of $\Theta_NQ\Theta_N$.

\subsection{Quantization of the torus}

We consider the standard torus, $X=\bbC/\Lambda$, $\Lambda = \bbZ^2$,
with the complex structure arising from that of $\bbC$.  A quantizing
line bundle on $X$ is holomorphically trivial when pulled-back to $\bbC$,
therefore its sections can be identified with entire functions on $\bbC$
satisfying a transformation law with respect to translations by elements
of $\Lambda$.  It is well-known that the functions that arise in this manner
are precisely theta functions. The quantizing line bundle is not unique,
since one can always tensor a given one with the flat line bundles over
$X$.  This gives rise to theta functions with characteristics.

\subsubsection{The Hilbert spaces}

A quantizing line bundle over $X$ can be constructed from a cocycle
$\chi: \bbC \times \Lambda\to \bbC\setminus\{ 0\}$ given by:
\begin{equation}\label{t5}
\chi(z\,,\,m+in)\,=\, (-1)^{mn}\,e^{\pi[z(m-in) +
\frac{1}{2}(m^2+n^2)]}
e^{-2\pi i[m\mu + n\nu]},
\end{equation}
where $\mu$ and $\nu$ are fixed real numbers (the so-called
characteristics of the bundle).  $\chi$ is called a cocycle
because it satisfies the condition
\begin{equation}\label{t4}
\chi(z,\lambda)\,\chi(z+\lambda, \mu)\,=\,\chi(z, \lambda + \mu).
\end{equation}
The quantizing line bundle is the quotient of $\bbC\times\bbC$ by
the equivalence relation:
\[
(z,a)\sim(w,b)\quad\Leftrightarrow\quad\exists \lambda\in\Lambda\
\mbox{such that} \ (w,b) = (z+\lambda\,,\, \chi(z,\lambda)a).
\]
For simplicity we will only consider here theta functions with
characteristics $(\mu,\nu)=(0,0)$.

We observe the following features of this construction:
\smallskip
\begin{enumerate}
\item The sections of this line bundle, $L$, are naturally identified with the
functions $f:\bbC\to \bbC$ such that
\begin{equation}\label{t6}
\forall (z,\lambda)\in\bbC\times\Lambda \qquad
f(z+\lambda)\,=\, \chi(z,\lambda)\,f(z).
\end{equation}
(Indeed the section associated to one such $f$ is defined by:
\[
s_f([z])\,=\,[(z,f(z))]
\]
where the square brackets denote equivalence classes.)
\item  For any integer $N$ the $N$-th power of $\chi$, $\chi^N$
is again a cocycle.  The line bundle it defines is the $N$-th
tensor power of $L$, $L^{\otimes N}$.
\item  A Hermitian structure on $L$ is defined by a function
$h : \bbC\to\bbR^{+}$ satisfying:
\begin{equation}\label{t7}
h(z)\,=\, |\chi(z,\lambda)|\, h(z+\lambda).
\end{equation}
The Hermitian metric we will consider is:
$ | [(z,a)]| = |a| e^{-\pi |z|^2/2}$.
\end{enumerate}

\begin{definition}
The space $\calH_N$ of holomorphic sections of the line bundle
$L^{\otimes N}$ is the space of entire functions $f: \bbC\to \bbC$ satisfying:
$\forall z\in\bbC\,,\ m+in\in\Lambda$
\[
f(z+m+in)\,=\,(-1)^{Nmn}\,e^{N\pi[z(m-in) + \frac{1}{2}(m^2+n^2)]}\,f(z).
\]
Its Hilbert space structure is given by the inner product
\begin{equation}\label{tb1}
<f\,,\,g>\,=\, \int_\calF\ f(z)\,\overline{g(z)}\,e^{-N\pi |z|^2}\,dx\,dy,
\end{equation}
where $\calF$ is a fundamental domain for $\Lambda$.
\end{definition}

The transformation law (\ref{t6}) is not the standard one for theta
functions (see \cite{Ba}, \cite{Mu}).  However, if $f\in\calH_N$, then
\[
f^x(z) := e^{-N\pi z^2/2}\,f(z)
\]
satisfies the classical transformation law
\[
f^x(z+m+in)\, =\,e^{N\pi(n^2 -2inz)}\,f^x(z).
\]
For future reference we also associate to $f\in\calH_N$ the function
\[
f^y(z) := e^{N\pi z^2/2}\,f(z)
\]
which satisfies the transformation law
\[
f^y(z+m+in)\, =\,e^{N\pi(m^2 +2mz)}\,f^y(z).
\]
Note that, in particular
\[
f^x(z+n) = f^x(z),\quad\text{and}\quad f^y(z+im) = f^y(z).
\]
Exploiting these periodicity conditions we now exhibit two (dual)
basis of $\calH_N$.

We begin with the functions $f^x$.  They can be expanded in Fourier
series,
\[
f^x(z) = \sum_{m=-\infty}^\infty\, a_m\,e(mz)
\]
where we let $e(z)\,:=\,e^{2\pi i z}$.  The transformation law for
$f^x$ becomes a relation among the Fourier coefficients, namely
$a_{m+Nn}\,=\,e^{-\pi(Nn^2+ 2mn)}\,a_m$.  This
shows that the dimension of the space of theta functions of order $N$
is $N$, as the values of $a_0,\ldots a_{N-1}$ determine the Fourier
series.  This leads to considering some special theta functions obtained
by letting exactly one of the coefficients $a_0,\ldots a_{N-1}$ be non-zero.
These functions give rise to a basis of $\calH_N$.  More precisely:
\begin{lemma}\label{Bax}
For $j=0,\ldots,N-1$, let $\vartheta_j^\N(z)$ be defined by the
Fourier series
\[
\vartheta_j^\N(z)\,=\,(2N)^{1/4}\, e^{N\pi z^2/2}
\sum_{n=-\infty}^\infty\,
e^{-\pi N(n+j/N)^2}\, e(z(j+Nn)).
\]
Then $\vartheta_j^\N\in\calH_N$, and the set
$\{\vartheta_j^\N,\ j=0\ldots,N-1\}$ is an orthonormal basis of $\calH_N$.
\end{lemma}

We can carry out a similar construction by considering the Fourier
series of the functions $f^y$ where $f\in\calH_N$.  This results in
a different basis of $\calH_N$:
\begin{lemma}\label{Bay}
For $k=0,\ldots,N-1$, let $\beta_k^\N(z)$ be defined by the
Fourier series
\[
\beta_k^\N(z)\,=\,(2N)^{1/4}\, e^{-N\pi z^2/2}
\sum_{m=-\infty}^\infty\,
e^{-\pi N(m+k/N)^2}\, e(-iz(k+Nm)).
\]
Then $\beta_k^\N\in\calH_N$, and the set
$\{\beta_k^\N,\ k=0\ldots,N-1\}$ is an orthonormal basis of $\calH_N$.
\end{lemma}

It's a beautiful fact that the matrix relating these two bases is the
discrete Fourier transform.  We claim that
\begin{equation}
\beta_k^\N(z) = \sum_{j=0}^{N-1}\, e^{-2\pi ikj/N}\,\vartheta_j^\N(z).
\end{equation}
We refer to \cite{Aus} and especially \cite{KP}, Proposition 3.17
for the reason for this.

\subsubsection{Matrix coefficients}
We now turn to a calculation of matrix coefficients of Toeplitz
operators on the torus.  We recall the following result of
\cite{BGPU} (Corollary 4.8):

\begin{lemma}
Let $(x,y)$ be standard coordinates on the torus, so that $z=x+iy$
in the previous formulas.  Let $f(x,y)$ be a symbol that is
in fact a smooth periodic function of $y$ alone.  Then the $\vartheta_j^\N$
are eigenvectors of the B-T operator $T_f^N$,
\[
T_f^N\vartheta_j^\N = \lambda_j^\N\vartheta_j^\N,\quad\text{where}\quad
\lambda_j^\N = \sum_{n=-\infty}^\infty a_n\,e^{-\pi n^2/2N}\,e^{-2\pi nj/N}
~ \sim f|_{y=-j/N}
\]
and where the $a_n$ are the Fourier coefficients of $f$ (with respect to
$y$).
\end{lemma}
Notice that in particular all such operators are normal.  Similarly,
the basis $\beta_k^\N$ consists of eigenvectors of any Toeplitz
operator with total symbol a function of $x$ alone.

Continuing our calculations of matrix coefficients,
let us now take a symbol $f$ of the form:
\begin{equation}\label{torus1}
f(x,y) = h(y)\,e^{2\pi i l x},\qquad\text{where}\quad l\in\bbZ.
\end{equation}

\begin{lemma}
If $f$ is of the form (\ref{torus1}), then
$\inner{f\vartheta_j^\N}{\vartheta_k^\N}$ is zero unless $k=(j+l)$ mod $N$,
in which case
\[
\inner{f\vartheta_j^\N}{\vartheta^\N_{(j+l)\,\text{\tiny mod}\,N}}
= h(-\frac{2j+l}{2N}) + O(1/N).
\]
More precisely, the matrix coefficient above is equal to
\[
e^{-\pi l^2/2N}\, (e^{-\pi\Delta_y/2N}\,h) (-\frac{2j+l}{2N})
\]
where $\Delta_y = \frac{d^2\ }{dy^2}$.
\end{lemma}
\begin{proof}
It is easy to verify the first statement, and, in case $k=j+l$ mod $N$,
one computes that
\[
\inner{f\vartheta_j^\N}{\vartheta^\N_{j+l\,\text{\tiny mod}\,N}} =
\sqrt{2N}\,\int_0^1\,\Psi_{j,l}(y)\,h(y)\,dy
\]
where
\[
\Psi_{j,l}(y) = e^{-\pi l^2/2N}\,\sum_{n=-\infty}^\infty\,
e^{{-2\pi N[y+n+(2j+l)/2N]}^2}.
\]
If $h=1$, then we have:
\[
\inner{e^{2\pi il x}\,\vartheta_j^\N}{\vartheta^\N_{(j+l)\,\text{\tiny mod
}\,N}}
= \sqrt{2N}\,e^{-\pi l^2/2N}\,
\sum_{n=-\infty}^\infty\,\int_0^1\, e^{{-2\pi N[y+n+(2j+l)/2N]}^2}\,dy
\]
\begin{equation}\label{torus.b}
=\sqrt{2N}\, e^{-\pi l^2/2N}\,\int_{-\infty}^\infty\,e^{-2\pi Ns^2}\,ds
= e^{-\pi l^2/2N}.
\end{equation}
This is $1+O(1/N)$ and therefore satisfies the desired estimate.
To proceed in case $h$ is not constant,
notice that by the Poisson summation formula $\Psi_{j,l}$ can be
written as
\[
\Psi_{j,l} = \frac{e^{-\pi l^2/2N}}{\sqrt{2N}}\,\sum_{k=-\infty}^\infty
\, e^{-\pi k^2/2N}\,e^{2\pi i k[y+\frac{2j+l}{2N}]}
\]
Therefore
\[
\inner{f\vartheta_j^\N}{\vartheta^\N_{(j+l)\,\text{\tiny mod}\,N}} =
e^{-\pi l^2/2N}\,\sum_{k=-\infty}^\infty\, e^{-\pi k^2/2N}\,
e^{\pi i k\frac{2j+l}{N}}\,\hat{h}(-k)\sim h(-\frac{2j+l}{2N}).
\]
\end{proof}
\begin{corollary}\label{CorT}
Let $k_1<k_2$ be two integers, and let $h_l(x)$ be smooth 1-periodic
functions, $k_1\leq l\leq k_2$.  Then there exists a B-T operator
on the torus, $\{T^\N\}$, with principal symbol
\[
f(x,y) = \sum_{l=k_1}^{l=k_2} h_l(y)\,e^{2\pi i l x}
\]
and such that for each $N$ the matrix entries of $T^\N$ in the basis
$\{\vartheta_j^\N\}$ are given by the formula
\begin{equation}\label{torus.c}
\inner{T^\N\vartheta_j^\N}{\vartheta_{(j+l)\,\text{\tiny mod}\,N}}  =
h_j(-\frac{2j+l}{2N}) + O(N^{-\infty}),
\end{equation}
where the estimate is uniform (in $j$, $l$).
\end{corollary}
\begin{proof}
For each $l$ there exists 
a sequence of periodic functions $f_{l,m},\ m=0,1,\ldots$ such that
\[
\Bigl(e^{-\pi l^2/2N}\, e^{-\pi\Delta_y/2N}\,\Bigr)
\sum_{m=0}^\infty N^{-m}\,f_{l,m}(y)
\sim h_l (y).
\]
where the left-hand side is considered a formal power series of $1/N$
(clearly $f_{l,0}(y) = h_l (y)$).
By the Borel summation method, there exists an $N$-dependent function
$f_l(y,N)$ periodic and smooth in $y$ such that
$f_l(y,N)\sim \sum_{m=0}^\infty N^{-m}\,f_{l,m}(y)$,
estimates the $C^\infty$ topology. 
By the previous Lemma, the matrix coefficients of the Toeplitz operator
with $N$-dependent multiplier
$f(x,y;N) = \sum_{l=k}^{l=m} f_l(y,N)\,e^{2\pi i l x}$
satisfy (\ref{torus.c}).  Since the sum over $l$ is finite the
estimates are uniform.
\end{proof}


\subsubsection{Localization of pseudomodes}
We now proceed to describe in concrete
terms the phase space localization of the pseudomodes constructed in \S 2.
We begin by pointing out that the $\vartheta^\N_j$ concentrate, as $N\to\infty$,
along the circles on the torus defined by $y=$constant.
This is because the $\vartheta_j$ are eigenvectors of Toeplitz operators
with symbols $f=f(y)$.  A precise statement is:

\begin{lemma} Fix $y_0=\frac{j_0}{N_0}$ a rational number modulo 1.
Then the microsupport of the sequence
$\{\vartheta^{(kN_0)}_{kj_0}\;;\;k=1,2,\ldots\}$
equals the circle $y=y_0$ on the torus.  The $\beta^\N_j$ accumulate
along circles $x=x_0$ in an analogous fashion.
\end{lemma}

An even more precise statement is that the sequences of this lemma are
Legendrian states associated to the corresponding circles, in the
sense of \cite{BPU}.

Consider now a pseudomode, $\psi_N$, associated with a Toeplitz
operator $T^\N_f$ on the
torus, with microsupport $(x_0,y_0)\in X$.  It follows from the previous
lemma that the components of $\psi_N$
\[
a^\N_j :=\inner{\psi_N}{\vartheta^\N_j}
\]
in the basis of the $\vartheta_j$ concentrate around the values $j$
where $j/N\approx y_0$.
Similarly, the coefficients
\[
b^\N_j :=\inner{\psi_N}{\beta^\N_j}
\]
concentrate around the values $j$ where $j/N\approx x_0$.
More precisely, the concentration occurs in a neighborhood of size
$O(\sqrt{N})$ of these values.  As mentioned, the sequence
$b^\N_j$ is the finite Fourier transform of $a^\N_j$.  Therefore,
{\em in these coordinates the pseudomode is localized on both sides
of the Fourier transform.}

\subsection{The complex projective line}
We begin by reviewing the quantization of $\CP$, for completeness.
The quantization of the complex projective line arises in connection
with the irreducible representations of SU$(2)$.  Recall that up to 
isomorphism such irreducible representations are those realized in the
spaces
\[
\calH_N := \{\, f(w_1,w_2)\;;\;  f\ \mbox{a homogeneous polynomial
of degree } N\,\}.
\]
Specifically, if $f\in\calH_N$ and $g\in$ SU$(2)$, then
\[
(g\cdot f)(w_1,w_2) = f(g^{-1}\cdot (w_1,w_2))
\]
where the action on the right-hand side is the natural action of
SU$(2)$ on $\bbC^2$.  Here $N=0,1,2,\cdots$ is a non-negative integer.
These representations are unitary if we put on $\calH_N$
the Hermitian inner product
\[
\inner{f_1}{f_2} = \int_{S^3}\,f_1\,\overline{f_2}\,dV_{S^3}\,,
\]
where $S^3\subset\bbC^2$ is the unit sphere and $dV_{S^3}$ is its
volume form.  We note without proof that the vectors
\begin{equation}\label{sph01}
\ket{j,N} = \sqrt{\frac{N+1}{\pi}}\,
\sqrt{C_N^j}\,w_1^{j}\,w_2^{N-j},\quad\ 0\leq j\leq N
\end{equation}
form an {\em orthonormal} basis of $\calH_N$ (consisting of eigenvectors
for the operator induced by $\sigma_3\in\text{su}(2)$, see below).

The circle group, $S^1\subset\bbC$, acts freely on $S^3$, by complex
multiplication.  Therefore $S^3$ is a circle bundle over the abstract
quotient,  $X := S^3/S^1$, which can be identified with the space
$\CP$ of all complex lines through the origin in $\bbC^2$.  There is
a natural Hermitian line bundle, $L^*$, over $\CP$ (the one whose
fiber at $\ell\in\CP$ is $\ell$ itself), and in fact $P=S^3\subset
L^*$ is the unit circle bundle.  The dual bundle, $L\to\CP$, is the
so-called hyperplane bundle and, as it turns out, quantizes $\CP$
(meaning that the natural connection on it has curvature the Fubini-Study
symplectic form of $\CP$). 
By homogeneity, the
functions on $\calH_N$ transform very simply along the orbits of $S^1$,
and the inner product above is invariant under $S^1$.  In fact we can 
regard the elements of $\calH_N$ precisely as the
holomorphic sections of $L^{\otimes N}$.  We are therefore exactly
in the setting of the previous sections.

The space $\CP$ is isomorphic to a two-dimensional sphere, 
as follows.  Define a map:
\begin{equation}\label{sph2a}
\Phi: \CP\to \mbox{su}(2)
\end{equation}
by the following rule:  For each $\ell\in\CP$, $\Phi(l)$
is the matrix having $\ell$ as an eigenspace, with associated 
eigenvalue $i/2$, and $\ell^\bot$ as another eigenspace, 
with associated eigenvalue $-i/2$.  (The choice of
spectrum is dictated by the normalization that the area
of $S^3/S^2$ agrees with the one induced by the Killing form.)  
One can show that, if we write $\ell = [w_1,w_2]\in \CP$,
where $(w_1,w_2)\in S^3$ is a representative, then the previous map is:
\begin{equation}\label{sph5}
\Phi([w_1,w_2]) = \frac{i}{2}
\begin{pmatrix}
|w_1|^2-|w_2|^2 & 2w_1\overline{w_2}\\
2w_2\overline{w_1} & |w_2|^2-|w_1|^2
\end{pmatrix}.
\end{equation}
It is a fact that $\Phi$ is a moment map for the natural SU$(2)$ action
on $\CP$.  Recall also that the matrices
\begin{equation}\label{sph4}
\sigma_1 = \frac{1}{2}
\begin{pmatrix}
0 & i \\
i & 0
\end{pmatrix}, \ 
\sigma_2 = \frac{1}{2}
\begin{pmatrix}
0 & -1 \\
1 & 0
\end{pmatrix}, \ 
\sigma_3 = \frac{1}{2}
\begin{pmatrix}
i & 0 \\
0 & -i
\end{pmatrix}
\end{equation}
form a standard orthogonal basis of su$(2)$ 
(such that $[\sigma_1,\sigma_2]=\sigma_3$, etc., and $\norm{\sigma_j}=1/2$)
if we give to su$(2)$ the \SU-invariant inner product
\[
\inner{A}{B} = -\frac{1}{2}\tr{AB}.
\]
Clearly $\Phi$ is an equivariant diffeomorphism onto its image, which,
geometrically, is the sphere of radius 1/2 and, algebraically, a
(co)adjoint orbit.  (It turns out that the symplectic form on $\CP$ is twice
the area form.  In general the symplectic form on an orbit of radius
$s$ is the area form divided by $s$, see \cite{Wood} pg.\ 54.) We will
henceforth identify $\CP$ with this sphere/co-adjoint orbit.

Let $x_j:\text{su}(2)\to\bbR$ be the $j$-th coordinate function,
$x_j(A) = \inner{A}{2\sigma_j}$, so that the image of $\Phi$ is the
sphere $\sum_{j=1}^3 x_j^2 = \frac{1}{4}$.
The description of $\CP$ as a
sphere means that we can speak of restrictions of linear functions 
from su$(2)$ to $\CP$.  In particular, a crucial role in what follows
will be played by the function $I:\CP\to [0,1]$ given by
\begin{equation}
I(\ell) = \inner{\Phi(\ell)}{2\sigma_3}+ \frac{1}{2}.
\end{equation}
The Hamilton flow of this function with respect to the natural symplectic
structure on $\CP$ (as a co-adjoint orbit) is given by the action of the
one-parameter subgroup $\exp(2t\sigma_3)$, and geometrically is rotation
around the $\sigma_3$ axis.  Accordingly, we introduce the polar angle,
$\theta$, regarded as a multivalued function on $\CP$ (undefined at the poles).
$(I,\theta)$ are action-angle coordinates on $\CP$, and, in particular, the
symplectic form on $\CP$ is 
\[
\omega = dI\wedge d\theta.
\]

Any Berezin-Toeplitz operator gives rise to a sequence of matrices,
namely, the matrices representing the operator in the canonical basis,
$\{\ket{j,N}\}$, of $\calH_N$.  
We will compute below (approximately) the matrix
of a Toeplitz operator with principal symbol a given function
$f:\CP\to\bbC$.  We begin with some remarks about such functions.
Let 
\begin{equation}\label{sph6}
f(I,\theta) = \sum_{l=-\infty}^\infty\, e^{il\theta}\,f_l(I)
\end{equation}
be the Fourier series expansion of $f$ with respect to the action
of $S^1$ by rotations around the $\sigma_3$ axis.
Although we are writing this expansion in in action-angle variables,
each summand is a smooth function on the sphere, which imposes
boundary conditions on the $f_l$.  Specifically, one can prove that 
the functions $f_l(I)$ such that $f_l(I) e^{il\theta}$ 
is smooth as a function on $\CP$ are of the form
\begin{equation}\label{sph7}
f_l(x) = 
(x(1-x))^{l/2}\, g_l(x) 
\end{equation}
where $g_l$ has a smooth extension to a neighborhood of $[0,1]$.
In fact, $x_1+ix_2 = (\frac{1}{4}-x_3^2)^{1/2}\,e^{i\theta}$, so
\[
(I(1-I))^{l/2}\,e^{il\theta} = (x_1+ix_2)^l.
\]
Therefore, if (\ref{sph7}) holds,
\[
f_l(I)e^{il\theta} = g_l(I)(x_1+ix_2)^l
\]
which is clearly smooth on the sphere if $g_l\in C^\infty[0,1]$.

\medskip
Let us now turn to the computation of matrix elements of B-T operators
on $\CP$.  By linearity of Toeplitz quantization, it will suffice to
compute the matrix of a Toeplitz operator with symbol
$e^{il\theta}\,f_l(I)$ for a given integer $l$.  We begin with the
case $l=0$.

\begin{lemma}\label{Sph2}
Let $\alpha$ be a smooth function on $[0,1]$.  Then there is a B-T
operator on the sphere, $\{A^\N\}$, 
with principal symbol $\alpha\circ I$, which
is diagonal in the standard basis of $\calH_N$ and whose 
$j$-th diagonal entry is $\alpha(j/(1+N))$, $0\leq j\leq N$.
\end{lemma}
\begin{proof} 
The sequence of operators, $Z = \{Z^\N:\calH_N\to\calH_N\}$ such that
\[
Z^\N \ket{j,N} = \frac{j}{N+1}\  \ket{j,N}, \qquad 0\leq j\leq N
\]
(where $\ket{j,N}$ is defined in (\ref{sph01})
is a Berezin-Toeplitz operator with symbol $I$
(see Lemma 3.4 and the ensuing discussion in \cite{BGPU}).
Let
\[
A^\N = \frac{1}{2\pi}\,\int\,e^{itZ^\N}\,\hat{\alpha}(t)\,dt,
\]
where $\hat{\alpha}$ is the Fourier transform of a compactly-supported
smooth extension of $\alpha$.  Since $Z$ is a self-adjoint Toeplitz 
operator of order zero, $\{e^{itZ^\N}\}$ is a unitary Toeplitz operator 
of order zero (see for example Proposition 12 of \cite{Ch}) and
symbol $e^{itx_3}$.  It is easy to check that $\{A^\N\}$ has the
desired properties.
\end{proof}

\medskip
To compute the matrix elements of operators with a smooth symbol 
$e^{il\theta}\,f_l(I)$ where $l\not=0$,  
we introduce the raising and lowering operators.
Let
\[
J_k^N := i\times \mbox{the operator induced by }\sigma_k\mbox{ on }\calH_N,\ 
k=1,2,3. 
\]
Then:
\[
J_1^N = \frac{1}{2}\Bigl( w_2\frac{\partial}{\partial w_1} + 
 w_1\frac{\partial}{\partial w_2}\Bigr)
\quad
J_2^N = \frac{i}{2}\Bigl( w_2\frac{\partial}{\partial w_1} - 
 w_1\frac{\partial}{\partial w_2}\Bigr)
\]
and
\[
J_3^N = \frac{1}{2}\Bigl( w_1\frac{\partial}{\partial w_1} - 
 w_2\frac{\partial}{\partial w_2}\Bigr).
\]
In particular, the lowering and raising operators,
$J^N_{\pm} = J_1^N\pm i J_2^N$, are:
\[
J^N_+ = w_1\frac{\partial}{\partial w_2},\qquad
J^N_{-} = w_2\frac{\partial}{\partial w_1}
\]
(and the vectors $\ket{j}$, $j=0,\ldots, N$
are eigenvectors of $J^N_3$, with eigenvalue $j-\frac{N-1}{2}$.)

\begin{lemma}
The matrix of $J_{-}^N$, resp.\ $J_{+}^N$,
in the standard basis has zero entries except along
the supra-diagonal, resp.\ infra-diagonal,
along which the entries are equal to
\[
m_j = \sqrt{j(N-j+1)}, \qquad j=1,\ldots, N.
\] 
Moreover, the sequence $\{\frac{1}{N}J_{\pm}^N\}$ is a B-T operator with
symbol $x_1\pm ix_2$.  
\end{lemma}
The first two statements follow a simple calculation; for the last
statement we refer to \cite{BGPU}.

We are now in a position to describe the matrices of the B-T
operators on the sphere:

\begin{proposition}
Let $f:\CP\to\bbC$ be a smooth function with a finite Fourier
series, (\ref{sph6}), in action-angle variables, and let
$f_l(x) = (x(1-x))^{l/2}\,g_l(x)$ with $g_l\in C^\infty[0,1]$. 
Let $\calM^\N_{\pm}$ be the matrix of $J_{\pm}^N$, described in the
previous Lemma.  Then there exists a B-T operator, $\{T^\N\}$, with symbol
$f$ and such 
that the sequence of matrices $\{\calT^\N\}$ of $\{T^\N\}$ in the
standard basis of $\calH_N$ satisfies:
\begin{equation}\label{sph6a}
\calT^\N = \sum_{l} \Bigl(\calM^\N_{\sgn(l)}\Bigr)^{|l|}
\,\calA^\N(g_l) + O(N^{-\infty})
\end{equation}
where $\calA^\N(g_l)$ is the $(N+1)\times (N+1)$ diagonal matrix
with diagonal entries $g_l(j/(N+1))$ and 
the estimate is in any matrix norm.
\end{proposition}
\begin{proof}
By linearity and the assumption that the sum (\ref{sph6}) is finite
it suffices to prove the proposition for $f$ of the form
$f(I,\theta) = e^{il\theta} f_l(I)$.  The case $l=0$ is covered by
Lemma \ref{Sph2}.  It suffices to consider the case $l>0$.
Applying Lemma \ref{Sph2}
to $g_l$ we obtain a diagonal B-T operator, $\{A^\N(g_l)\}$, 
with diagonal entries $g_l ((j-1)/N)$, $0\leq j\leq N$. 
It is clear that the B-T operator 
\[
T^\N := \Bigl(\frac{1}{N}J_-^\N\Bigr)^l\circ A^\N
\]
has the desired properties.
Notice that one can choose any order of the products appearing
in (\ref{sph6a}) and still find a B-T operator with the
desired properties.
\end{proof}

\subsubsection{Linear Hamiltonians}
Recall the moment map, $\Phi:\CP\to\mbox{su}(2)$ given by 
(\ref{sph5}).  Given $M\in\mbox{sl}(2,\bbC)$, we can pull-back
by $\Phi$ the complex-valued linear function on su$(2)$,
\[
\mbox{su}(2)\ni A\mapsto -\frac{1}{2}\tr (AM).
\]
Let us denote the pull-back by $F_M:\CP\to\bbC$; specifically
\[
\forall\ell\in\CP\qquad F_M(\ell) = -\frac{1}{2}\tr (\Phi(\ell)M).
\]
We will call functions such as $F_M$ {\em linear Hamiltonians}.
Since $\Phi$ is a moment map, the assignment $M\to F_M$ is a Lie
algebra morphism:
\[
\forall M_1,\ M_2\in \mbox{sl}(2,\bbC)\qquad F_{[M_1,M_2]} =
\PB{F_{M_1}}{F_{M_2}}.
\]
In particular, if we continue to denote by
$x_j:\CP\to\bbR$ the restriction to $\CP$ of the coordinate 
functions, then we have the identity:
$\PB{x_1}{x_2}=x_3$, and also its cyclic permutations.

If $M$ is semi-simple, there exists $g\in\mbox{SU}(2)$ such that
$gMg^{-1}$ is diagonal, i.e.
\[
\exists g\in\mbox{SU}(2),\ \mu\in\bbC\setminus{0}\qquad
gMg^{-1} = \mu\,\sigma_3.
\]
Here's a very concrete example.  Take
\[
V=\exp(it\sigma_2) =
\begin{pmatrix}
\cosh(t/2) & -i\sinh(t/2)\\
i\sinh(t/2) & \cosh(t/2)
\end{pmatrix}
\]
and consider $A= A(t)\in\text{sl}(2,\bbC)$ equal to
\[
A = V\sigma_1 V^{-1} = i\sinh(t)\,\sigma_1 +
\cosh(t)\,\sigma_3.
\]
The classical Hamiltonian, $F_A:\CP\to\bbC$ is $F_A =
i\sinh(t)\,x_1+\cosh(t)\,x_3$, and its image is the 
interior of an ellipse,
\begin{equation}\label{exa}
\frac{x^2}{\cosh(t)^2} + \frac{y^2}{\sinh(t)^2} \leq \frac{1}{4}.
\end{equation}
If we let $T^\N:\calH_N\to\calH_N$ be $1/N$ times the operator
image of $\frac{1}{i}A$ by the representation $\rho_N$, then,
for all $t$, $T^\N$ is diagonalizable with real spectrum 
$\{\frac{j}{N}-\frac{N-1}{2N}\;;\;j=0,\ldots,N-1\,\}$, which is
contained in the major axis of the image of $F_A$.

Notice that $\PB{\Re F_A}{\Im F_A} = \cosh(t)\sinh(t)\,x_2$,
and therefore for every point in the image of $F_A$ there is exactly
one $x$ where this Poisson bracket is negative.  Therefore,
for each $\lambda$ in the interior of the image the norm of the
resolvent $\norm{(T^\N-\lambda I)^{-1}}$ is $O(N^\infty)$.

The level sets of the norm of the resolvent for this example
resemble the equipotential curves of a uniform electric charge
distribution on the line segment $[-\frac{1}{2},\frac{1}{2}]$.
Notice, incidentally, that the image under $F_A$ of the
level curves of the Poisson bracket $\PB{\Re F_A}{\Im F_A}$
are ellipses crossing this line segment. 
{\em This example therefore
shows that there is not a direct relationship between the norm
of the resolvent and the size of the Poisson bracket of the
real and imaginary parts of the symbol.}

Finally, notice that every $\lambda$ on the boundary of the elliptical
region (\ref{exa}) satisfies the hypotheses of Theorem \ref{II}.  The
set $f^{-1}(\lambda)$ consists of exactly one point, where $\PB{\Re
F_A}{\Im F_B} = \cosh(t)\sinh(t) x_2=0$.  However, one of the double
brackets involving real and imaginary parts of $F_A$ is non-zero at
this point.  Thus Theorem \ref{II} applies, with $k=2$.  Notice that,
in spite of the estimates (\ref{1i}) the distance from $\lambda$ to
spectrum is $O(1)$.  This is in agreement with Theorem 3 of \cite{DSZ}:
Thinking of $\CP$ as a {\em real} manifold $\calO$, the symbol $f$ has
an obvious holomorphic extension to a complexification of $\calO$,
namely, the (co)adjoint orbit of SL$(2,\bbC)$ through $\sigma_3$.

\section{Final Remarks}

\subsection{On the numerical range}

Let $T_f=\{T^\N\}$ be a Berezin-Toeplitz operator with symbol $f:X\to\bbC$.
For each $N$, the numerical range of $T^\N$ is the set
\begin{equation}
W_N:=\{\,\inner{T^\N\psi}{\psi}\;;\;\psi\in\calH_N\,,\,\norm{\psi} = 1\,\}.
\end{equation}
We define $W_\infty$ as the limit of the ranges $W_N$ as $N\to\infty$:

\begin{definition}
A complex number $\lambda$ is in $W_\infty$ iff for all $\epsilon >0$
there exists $K>0$ such that for all $N>K$
\[
\Delta_\epsilon(\lambda)\cap W_N\not=\emptyset,
\]
where $\Delta_\epsilon(\lambda)$ is the disc of radius $\epsilon$ centered
at $\lambda$.
\end{definition}

\begin{proposition}
$W_\infty$ is the convex hull of the image of the classical symbol, $f$.
\end{proposition}
\begin{proof}
It is well-known that, for each $N$, $W_N$ is convex (see \cite{GR}).
It follows easily that $W_\infty$ is convex as well.  Moreover, if
$x\in X$ and $\psi^N_x$ is a coherent state at $x$, then
\[
\frac{\inner{T^\N\psi^N_x}{\psi_x^N}}{\inner{\psi^N_x}{\psi^N_x}}\to f(x).
\]
This shows that $W_\infty$ contains the image of $f$.

To show that $W_\infty$ is actually the convex hull, consider a line
of equation $ax+by=c$, and let $\lambda\in W_\infty$. 
By definition, there exists a sequence $\lambda_N\in W_N$
converging to $\lambda$, and therefore there exists a sequence of
unit vectors $\{\psi_N\in\calH_N\}$ such that 
\[
\lambda_N = \inner{T^\N\psi_N}{\psi_N}\to\lambda.
\]
Let us write $\lambda_N = x_N + iy_N$ for the real and imaginary
parts of $\lambda_N$, and $f=p_1+ip_2$ for the real and imaginary parts
of the symbol $f$.  Then
\[
ax_N + by_N = \int_X (ap_1 + bp_2)\,|\psi_N|_s^2\,dm + O(1/N).
\]
Assume that the region: $ax+by>c$ does not intersect the image of $f$.
Then for all $\in X$ $a\phi(m)+bg(m) \leq c$, and
therefore
\[
ax_N+by_N \leq c\,\int_X |\psi_N|^2\,dm_N + O(1/N)  = c + O(1/N).
\]
Letting $N\to\infty$, we obtain that
$ax_\infty+by_\infty\leq c$, where $\lambda = x_\infty+iy_\infty$.
Thus $\lambda$, and therefore all of $W_\infty$,
is on the same side of the line $ax+by=c$ as the image of $f$.
\end{proof}

\subsection{The weak Szeg\"o limit theorem}  

In the non-selfadjoint case, one has the following version of the
Szeg\"o limit theorem for B-T operators:

\begin{proposition}\label{4I}
Let $T_f=\{T^\N\}$ be a B-T operator with principal symbol
$f:X\to\bbC$, and let $F(z)$ be a function of a complex variable
analytic on a simply-connected region containing the image of $f$. 
Then
\begin{equation}\label{4a}
\frac{1}{\dim \calH_N} \tr F(T^\N) = \frac{1}{\text{Vol }X}\,
\int_X F\circ f\, dm + O(1/N),
\end{equation}
where $dm$ is the Liouville measure of $X$.
\end{proposition}
\begin{proof}
The idea of the proof is standard; we include some details
for completeness and to verify that the usual proof
is valid in the current setting.
As in the theory of pseudodifferential operators,
for each $\lambda$ not in the image of $f$ one can construct a
B-T operator, $B_\lambda$, such that the Schwartz kernel
of 
\[
R^\N_\lambda:= B^\N_\lambda\circ (T^\N-\lambda) - I_{\calH_N}
\]
is a smooth section of the bundle $\text{Hom}(L^N,L^N)\to X$ which
is rapidly decreasing in $N$ (together with all its derivatives).
The principal symbol of $\{B^\N_\lambda\}$ is $(f-\lambda)^{-1}$.
Multiplying on the right by $(T^\N-\lambda)^{-1}$ we obtain
\[
(T^\N-\lambda)^{-1} = B_\lambda + S^\N_\lambda
\]
where $\{S^\N\}$ has the same properties as $R^\N_\lambda$.
Let $\Gamma$ be a simple closed curve, positively oriented, contained
in a region where $F(z)$ is analytic and containing the image of $f$.
Then
\[
F(T^\N) = \frac{1}{2\pi i}\,\oint_\Gamma F(\lambda)\,B^\N_\lambda\,d\lambda
+\frac{1}{2\pi i}\, \oint_\Gamma F(\lambda)\,S^\N_\lambda\,d\lambda .
\]
The trace of the second term on the right-hand side is $O(N^{-\infty})$,
while 
\[
\frac{1}{\dim \calH_N} \tr\oint_\Gamma F(\lambda)\,B^\N_\lambda\,d\lambda
=\frac{1}{\dim \calH_N} \oint_\Gamma F(\lambda)
\,\tr\,B^\N_\lambda\,d\lambda.
\]
But it is known that $\frac{1}{\dim \calH_N}\tr\,B^\N_\lambda = 
\frac{1}{\text{Vol }X}\,\int_X (f-\lambda)^{-1}\, dm + O(1/N)$,
where the estimate is uniform for $\lambda$ on compact sets away from 
the image of $f$.
\end{proof}

\bigskip
In particular, if $\lambda$ is a complex number away from the image
of $f$, one has:
\begin{equation}\label{add9}
\frac{1}{\dim \calH_N} \tr [(T^\N-\lambda)^{-1}] = 
\frac{1}{\text{Vol }X}\,\int_X \frac{1}{f-\lambda}\,dm + O(1/N).
\end{equation}
Clearly the left-hand side of this equation is a sequence of
analytic functions in $\lambda$ defined away from the union of the
spectra of the $T^\N$.  On the other hand, the integral
$\int_X \frac{1}{f-\lambda}\,dm$ is analytic away from the image of $f$.
Examples show that the spectral radius of $T^\N$ has a limit, $R$, such
that the image of $f$ is {\em not} contained in the circle of radius $R$. 
It is not immediate to extend (\ref{add9}) to $\lambda$ with $|\lambda|>R$
but inside the image of $f$.

\appendix

\section{Hermite distributions and symbol calculus}

\noindent{\em Oscillatory integrals and symbols}

We place ourselves in the setting of \S 2.2:  Let 
$\calR\subset\calZ$ be a closed conic isotropic submanifold, and
$Y\hookrightarrow X$ the reduced isotropic submanifold of $X$.
Hermite distributions in $I^m(P,\calR)$ are defined locally as  
oscillatory integrals.  
To write down an explicit form for these integrals, we'll choose Darboux
coordinates $(q,p)\in \bbR^{2n}$ for $X$ with $q = (q',q'') \in \bbR^l\times
\bbR^{n-l}$ such that $Y = \{q'' = p = 0\}$. Here $l = \dim Y$.
For $P$ we then have coordinates
$z = (q,p,\theta)$, and we can always find a function $h(p,q)$ such that 
$\alpha|_{p=0} = d\theta - dh|_{p=0}$.  The lift of $Y$ to the isotropic 
$\calR$ is given by specifying that $\theta = h(q',0)$.  

Now we'll introduce phase coordinates $(\tau, \eta_1, \eta_2) \in \bbR_+ 
\times \bbR^{n-l} \times \bbR^n$, and a phase function
$$
\phi(z,\tau,\eta) = \tau(\theta - h(q,p)) + \eta_1\cdot q'' + \eta_2 \cdot p, 
$$
which parametrizes $\calR$.  A distribution $u\in I^m(P,\calR)$ can be written locally
as
$$
u(z) = \int e^{i\phi(z,\tau,\eta)} a(x,\tau,\eta/\sqrt{\tau}) \>d\tau\>d\eta,
$$
where the amplitude $a(z,\tau, u)$ is rapidly decreasing in $u$ and has a expansion in
$\tau$ of the form 
$$
a(z,\tau,\eta) \sim \tau^{m-1/2} \sum_{j=0}^\infty \tau^{-j/2} a_j(z,\eta).
$$

The symplectic spinor symbol $\sigma(u)$ should be thought of as $a_0$ written in a 
suitably invariant way.  The choice of $\phi$ defines for each $\rho\in\calR$
a canonical isomorphism between
$T^*(\bbR^{n-l} \times \bbR^n)$ and the symplectic normal 
$N_\rho$, by which the $H_\infty(N_\rho)$ portion of the symbol may
be pulled back to a rapidly decreasing function of $\eta$.  This gives the
$\eta$ dependence of $a_0$, while the half-form portion of the symbol encodes
the dependence on $z\in P$ in a coordinate-independent way.

The construction of $u$ in \S\ref{Pf} gives $u$ such that $\sigma_u = \nu_u
\otimes \kappa_u \otimes e$.  Here $e$ is a Gaussian in
$H_\infty(Z_\rho^\circ)$ which can be explicitly computed in terms of the
metric.  The component $\kappa_u \in H_\infty(E_\rho)$ satisfies 
$\calL(\kappa_u) = 0$ according to the construction.  This does
not fix $\kappa_u$, but we are
free to assume that $\kappa_u$ is also Gaussian.
Here $\calL$ is the operator on $H_\infty(E_\rho)$ given by the Heisenberg 
representation of the Hamiltonian vector field $\xi$ of $\sigma_Q$.  
But according to Theorem 11.4 of \cite{BG}, $\xi_\calZ$ 
is just the lift of the Hamiltonian vector field $\Xi_f$ on $X$
up to $\calZ$.   So $\calL$ could be written explicitly in terms of $f$.

By this construction of $u$, the leading amplitude
$a_0(z,\eta)$ is a Gaussian function of $\eta$:
$$
a_0(z,\eta) = c(z) e^{-\eta^t A^{-1} \eta/2}
$$
where $A$ is a symmetric matrix with positive definite real part.
Arguing as in the proof of Theorem 3.12 in \cite{BPU}, we can write
$u_N = \Pi_N(u)$ as 
$$
u_N((q,p,h(q,p)) =
e^{iNh(q,p)} \int e^{-iN\theta} e^{i\phi} a(z,\tau,\eta/\sqrt{\tau})
\;d\theta\> d\tau \>d\eta,
$$
(near $q'' = p = 0$) and apply stationary phase to the $N\to\infty$ limit.
The asymptotic result is that 
$$
u_N(q,p,h(q,p)) \sim c(q') N^{n-l/2-1/2} e^{-N(q'',p)A(q'',p)^t/2},
$$
where $c(q')$ is independent of $N$.
In particular, at a point $z_0\in P$ lying above $Y$, we have 
\begin{equation}\label{psirho}
u_N(z_0) \sim c(z_0) N^{n-l/2-1/2}.
\end{equation}

\bigskip
\noindent{\em A review of the Hermite calculus}

The composition of Hermite distributions is described in Theorem 9.4
of \cite{BG}, which implies that
$$
\Pi: I^m(R,\calR) \to I^m(R,\calR).
$$
The symbol calculus corresponding to this composition is based on symplectic
linear algebra found in \S6 of \cite{BG}, where the reader can find full details.  
Since the symbol calculus is somewhat involved, 
we begin with a review of the general details of the symbol map before applying 
it to our case.  

Let $V$ and $W$ be symplectic vector spaces, $\Gamma\subset V\times
W^{-}$ a Lagrangian subspace and $\Sigma\subset W$ an isotropic
subspace.  We think of $\Gamma$ as a canonical relation from $W$ to
$V$; $\Gamma\circ\Sigma$ is an isotropic subspace of $V$.  We will make
the simplifying assumption that
\begin{equation}\label{calc1}
U_0\,=\,\{\,w\in\Sigma\,;\,(0,w)\in\Gamma\,\} = 0,
\end{equation}
valid in the applications of the calculus to this paper.

We assume given a symplectic spinor on $\Sigma$ and a half-form on
$\Gamma$.  Recall that if $H_\infty(V)$ denotes the space of $C^\infty$
vectors in the metaplectic representation of the metaplectic group of
the symplectic vector space $V$, the space of symplectic spinors on
$\Sigma$ is $H_\infty(\Sigma^{\circ}/\Sigma)\otimes\hffm (\Sigma)$. 

Under the assumption (\ref{calc1}),
the (linear) symbol map of the Hermite calculus is a linear map

\begin{equation}\label{hmap}
H_\infty(\Sigma^{\circ}/\Sigma)\otimes \hffm (\Sigma) \otimes \hffm (\Gamma)
\to
H_\infty( (\Gamma\circ\Sigma)^{\circ}/ \Gamma\circ\Sigma)
\otimes \hffm ( \Gamma\circ\Sigma ).
\end{equation}

Our first goal here is to describe the map (\ref{hmap}). There are two
ingredients in its construction, which will be examined separately.
First however we must introduce the following vector spaces:

\begin{definition}\label{Udef}
$\displaystyle{
U_1\,:=\,\{\,w\in\Sigma^{\circ}\,;\,(0,w)\in\Gamma\,\}\,\subset W\,}$,
and
\[
U\,:=\,\mbox{image of}\ U_1\ \mbox{in}\ \Sigma^\circ/\Sigma\,\cong U_1.
\]
\end{definition}
These spaces enter the calculus in the following way:
\begin{lemma}\label{Calc1}  The subspace
\[
U\subset \Sigma^{\circ}/\Sigma
\]
is isotropic, and there is a natural identification
\begin{equation}
U^{\circ}/U \cong (\Gamma\circ\Sigma)^{\circ}/ \Gamma\circ\Sigma\,.
\end{equation}
\end{lemma}

\medskip
The first ingredient in the symbol map is a canonical isomorphism:
\begin{lemma}  Under (\ref{calc1}), there exists a canonical isomorphism
\[
\hffm (\Sigma)\otimes \hffm (\Gamma)\;\cong\;
\hffm (U_1^{\circ}) \otimes  \hffm ( \Gamma\circ\Sigma )\,.
\]
\end{lemma}
\begin{proof}
Let
\begin{equation}
\rho : \Gamma \oplus \Sigma \to U_1^{\circ}
\end{equation}
be the map $\rho ((v,w),w_1)\,=\,w-w_1$.
One can show that the image of this map is exactly $U_1^{\circ}$. 
Moreover, because of (\ref{calc1}), the projection
$((v,w),w_1)\mapsto v$ is an isomorphism
\[
\ker(\rho)  \cong \Gamma\circ\Sigma.  
\]
This is the non-trivial vertical arrow in the diagram:
\begin{equation}
\begin{array}{ccccc}
 & 0 & & & \\
 & \downarrow & & & \\
0  \to & \ker (\rho) & \to \Gamma \oplus \Sigma \to & U_1^{\circ}\to & 0 \\
 & \downarrow & & & \\
 & \Gamma\circ\Sigma & & & \\
 & \downarrow & & & \\
 & 0 & & & 
\end{array}\label{rhoseq}
\end{equation}
The horizontal sequence is just the natural short exact sequence associated
to the surjection $\rho$.  

Having established the existence of these exact sequences the desired
isomorphism follows from the behavior of the functor $\hffm$ 
when applied to short exact sequences and to direct sums.
\end{proof}

\medskip
The second ingredient in the Hermite calculus
is the following:
\begin{lemma}  Under the assumption (\ref{calc1}),
there exists a canonical map
\[
\calS (\Sigma^{\circ}/\Sigma)\to \nchffm (U_1)\otimes
\calS ( (\Gamma\circ\Sigma)^{\circ}/
\Gamma\circ\Sigma)\,.
\]
\end{lemma}
\begin{proof}

This is based on Lemma (\ref{Calc1}) and the following
generalization of a map defined by Kostant:

\medskip
\noindent {\bf Claim:}
Let $A$ (in our case we will take $A =\Sigma^{\circ}/\Sigma$)
be a symplectic vector space,
and $U\subset A$ an isotropic subspace.  Then there is a natural 
map
\begin{equation}
H_\infty (A) \to \nchffm (U)\otimes H_\infty (U^{\circ}/U)\,.
\end{equation}

\medskip\noindent
The desired map follows from these two claims, if we recall that
$U=U_1$ (because of (\ref{calc1}).
\end{proof}

To obtain the symbol map (\ref{hmap}), tensor the maps from the lemmas and use
the fact that the symplectic form on $W$ defines a natural identification
\[
\nchffm (U_1)\otimes \hffm (U^{\circ}_1)\cong \bbC\,.
\]

\bigskip\noindent{\em Proof of Proposition \ref{WpI}.}

Given $v \in I^m(P,\calR)$, we want to calculate the symbol of $u = \Pi(v)
\in I_\Pi^m(R,\calR)$.
In order to apply the symbol calculus reviewed above, we need to rewrite $\Pi(v)$
as the composition of a Lagrangian distribution with a Hermite.  Thus we introduce
$\pi: P \times P \to P$, the projection through the left factor, and 
$F:P\times P \to P\times P \times P$ the map $F(p_1,p_2) = (p_1,p_2,p_2)$.
We can then write
$$
\Pi(v) = \pi_* F^* (\Pi \boxtimes u),
$$
where (abusing notation slightly) $\Pi$ here denotes the integral kernel.

For the symbol calculation it suffices to localize to $\rho \in \calR 
\subset T^*P$. 
To simplify notation we will introduce vector spaces
$$
V = T_\rho(T^*P), \qquad W = V\times V\times V.
$$
In $W$ two copies of $V$ carry the opposite symplectic form, but to simplify notation
we just denote the vector space.  We define the vector spaces $R_\rho$, $Z_\rho$,
$E_\rho$, etc. as in \S\ref{Pf}.
Note that $Z_\rho$ is a symplectic subspace of $V$, and $R_\rho$ is isotropic in $Z_\rho$.
Locally, the symbol of $v \in I(P,\bbR)$ can be written
\begin{equation*}\begin{split}
\sigma(v)|_\rho &= \nu \otimes \kappa \otimes \lambda  \\
&\in  \hffm(R_\rho) \otimes H_\infty(E_\rho) \otimes H_\infty(Z_\rho^\circ).
\end{split}\end{equation*}
As an integral kernel,
$\Pi \in I(P\times P, \calZ_\rho^\Delta)$ with symbol
$$
\sigma(\Pi) = \sqrt{dz} \otimes e\otimes \bar e \in \hffm(Z_\rho) \otimes
H_\infty(Z_\rho^\circ) \otimes H_\infty(Z_\rho^\circ),
$$
where $dz$ is the canonical volume form given by the symplectic form on $Z_\rho$.

The operator $\pi_* F^*$ is a Lagrangian FIO with canonical relation
$$
\Gamma = \{(v;v,w,w);\; v,w \in V\} \subset V\times W.
$$
(To see that it's Lagrangian one needs to
keep track of the signs of the symplectic forms).
The combination $\Pi \boxtimes u$ is a Hermite distribution associated to the isotropic
$$
\Sigma = \{(z,z,y);\; z\in Z_\rho, y\in R_\rho\}.
$$
Note that $\Gamma \circ \Sigma = R_\rho$.

First we'll describe the spaces that play a role in Lemmas A1--4.  To begin, note that $R_\rho^\circ$
was calculated in Lemma 2.3 to be $Z_\rho^\circ \oplus E_\rho$, where $E_\rho$ is the symplectic normal of
$R_\rho$ as a subspace of $Z_\rho$.  Then we see that
$$
\Sigma^\circ/\Sigma = Z_\rho^\circ \times Z_\rho^\circ \times (Z_\rho^\circ \oplus E_\rho).
$$
Then
$$
U_1 = \{(0,v,v);\; v\in Z_\rho^\circ\} \subset W,
$$
and $U$ is the same set as a subspace of $\Sigma^\circ/\Sigma$,
so that $U^\circ/U \cong Z_\rho^\circ \oplus E_\rho = R_\rho^\circ/R_\rho$.  

Introducing the map $\rho$ as above, we see that
$$
\operatorname{Image}(\rho) = U_1^\circ = U_1 \oplus (V \times Z_\rho \times Z_\rho) \subset W,
$$
and $\ker \rho \cong R_\rho$.

With these identifications, we can decompose the pieces in the symbol map (A.2).
First of all, in Lemma A.3 the map reduces to the identity map on
$$
\hffm(Z_\rho) \otimes \hffm(R_\rho) \otimes \hffm(V), 
$$
combined with the obvious decomposition
$$
\hffm(V) \cong \hffm(Z_\rho) \otimes \hffm(Z_\rho^\circ).
$$
Also note that since $\Sigma^\circ/\Sigma = Z_\rho \times Z_\rho \times (Z_\rho\oplus E_\rho)$ we have
$$
H_\infty(\Sigma^\circ/\Sigma) = H_\infty(Z_\rho^\circ) \otimes H_\infty(Z_\rho^\circ) \otimes 
H_\infty(Z_\rho^\circ) \otimes H_\infty(E_\rho).
$$

This means we can break the symbol map (A.2) into three pieces.
The first is the identity map  
$$
H_\infty(E_\rho) \otimes \hffm R_\rho \to H_\infty(E_\rho) \otimes \hffm R_\rho.
$$
The second is a canonical pairing (A.6) applied to the symplectic space
$Z_\rho^\circ\times Z_\rho^\circ$ with the diagonal as Lagrangian subspace:
$$
H_\infty(Z_\rho^\circ) \otimes H_\infty(Z_\rho^\circ) \cong \nchffm(Z_\rho^\circ).
$$
Finally, the third is the natural identification
$$
\nchffm(Z_\rho^\circ) \otimes \hffm(Z_\rho^\circ) \otimes \hffm(Z_\rho) \otimes \hffm(V) \cong
\bbC,
$$
defined by the symplectic forms.

Applying this to the symbol of $\Pi(u)$, the third piece shows some natural
half-forms canceling, the second gives contributes a factor
$\langle e, \lambda\rangle$
which can be absorbed by changing
the half-form component $\nu\mapsto \nu'$.  And the first map then
gives the stated conclusion:
$$
\sigma(\Pi(v)) = \nu' \otimes \kappa \otimes e.
$$

\end{document}